\documentclass[12pt,reqno]{amsart}

\usepackage{amsmath}\allowdisplaybreaks  
\usepackage{amssymb}        
\usepackage{amsfonts}       
\usepackage{amsthm}         
\usepackage{mathrsfs}       
\usepackage{cases}          
\usepackage{siunitx}        
\usepackage{natbib}
\usepackage{array}
\usepackage{booktabs}       
\usepackage{longtable}      
\usepackage{multirow}
\usepackage{makecell}       
\usepackage{tabularx}       
\usepackage{nicematrix}     

\usepackage{tikz}\usetikzlibrary{tikzmark, arrows.meta}   
\usepackage{float} 
\usepackage{xcolor}           
\usepackage{graphicx}         

\usepackage[linesnumbered, ruled, vlined]{algorithm2e}  

\usepackage[colorlinks=true, linkcolor=blue, citecolor=blue, urlcolor=blue]{hyperref} 
\usepackage[nameinlink]{cleveref}	
\usepackage{color}
\usepackage{datetime}
\usepackage[left=2cm, right=2cm, top=2cm, bottom=2cm]{geometry}
\usepackage{cleveref}
\usepackage{ulem}           

\makeatletter
\newcommand\figcaption{\def\@captype{figure}\caption} 
\newcommand\tabcaption{\def\@captype{table}\caption}
\makeatother

\newtheorem{thm}{Theorem}[section]

\newtheorem{example}{Example}[section]

\numberwithin{equation}{section}
\usepackage{natbib}
\setcitestyle{numbers,square,comma}
\begin{document}
	\title[High Accuracy Techniques Based AFEM]{High Accuracy Techniques Based Adaptive Finite Element Methods for Elliptic PDEs}  
	\author[J. Xiao, Y. Liu and N. Yi]{Jingjing Xiao $^a$, Ying Liu $^{b, *}$ and Nianyu Yi$^{c}$}
	\thanks{$^*$ Corresponding author.}
	\thanks{$^a$ School of Mathematics, Shandong University, Jinan, 250100, P.R.China, Email: jingjing\_xiao@mail.sdu.edu.cn}
	\thanks{$^b$ Corresponding author. School of Science, Xi’an University of Technology, Xi'an 710048, P.R.China, Email: liuyinglixueyuan@xaut.edu.cn}
	\thanks{$^c$ School of Mathematics and Computational Science, Xiangtan University; Hunan National Center for Applied Mathematics; Hunan Key Laboratory for Computation and Simulation in Science and Engineering, Xiangtan 411105, P.R.China, Email: yinianyu@xtu.edu.cn.}
	
	\begin{abstract}
		This paper aims to develop an efficient adaptive finite element method for the second-order elliptic problem. 
		Although the theory for adaptive finite element methods based on residual-type a posteriori error estimator and bisection refinement has been well established, in practical computations, the use of non-asymptotic exact of error estimator and the excessive number of adaptive iteration steps often lead to inefficiency of the adaptive algorithm.
		We propose an efficient adaptive finite element method based on high-accuracy techniques including the superconvergence recovery technique and high-quality mesh optimization. The centroidal Voronoi Delaunay triangulation mesh optimization is embedded in the mesh adaption to provide high-quality mesh, and then assure that the superconvergence property of the recovered gradient and the asymptotical exactness of the error estimator. A tailored adaptive strategy, which could generate high-quality meshes with a target number of vertices, is developed to ensure the adaptive computation process terminated within $7$ steps. The effectiveness and robustness of the adaptive algorithm is numerically demonstrated.
	\end{abstract}
	
	
	\maketitle
	
	\section{Introduction}
	The adaptive finite element method (AFEM) \cite{BCNV2024} has been extensively studied and become a powerful tool in scientific computation, particularly effective for problems with singularities or multi-scale characteristics. 
	AFEM aims to get a numerical solution within a prescribed tolerance by using an optimal number of degrees
	of freedom. 
	Starting with an initial coarse mesh, the standard adaptive finite element algorithm creates a sequence of adapted meshes, and the corresponding finite element solutions are computed on each mesh. 
	
	An essential ingredient of AFEM is the error estimator, which provides global and local information on the numerical error. Globally, the a posteriori error estimator can be used as a stopping criterion to determine whether the finite element solution is an acceptable approximation. Locally, it can be used as an indicator to show the error distribution and guide local mesh adaption. Therefore, one can calculate the estimators to find out where the approximation error is large. Then instead of globally refining every element, one can select a subset of elements for which the error estimator is relatively large and refine these elements. 
	Several types of a posteriori estimators have been proposed for different problems (see, e.g., \cite{br, Babu_1978, Bernardi2000, BCNV2024, cz, cn, yan2001, Karel2010, fa, LXYC2024, yz, zz} and the references cited therein). 
	
	Another essential ingredient of AFEM is the mesh refinement method. 
	Mesh adaptivity techniques are generally divided into three categories: adaptive re-meshing methods \cite{Huang2011, Ju2006, DJ1990, PJ1992}, element subdivision methods \cite{D1996, AD2000}, and fixed-order mesh modification methods \cite{babuska2012, joe1989}. Adaptive re-meshing methods create an adaptive mesh using automatic mesh generation algorithms according to specified element size and shape. Element subdivision methods, such as regular refinement algorithm \cite{bank1998}, and bisection refinement algorithm \cite{rivara1984}, can preserve the shape regularity of elements. The third mesh adaption technique applies local mesh modifications in a fixed order, where the mesh quality \cite{babuska2012, joe1989} is improved by exploiting local mesh modification operations: swap, collapse, split, and relocation. 
	
	For the standard AFEM with residual type a posteriori error estimator and the bisection refinement, the convergence and optimality of the adaptive algorithm are well established \cite{BCNV2024, D1996, MKN2005, mn}. While these methods have been shown to be very successful in numerical solving partial differential equations, the computing efficiency
	could be poor performance in practice, manifested in the inaccuracy of error estimators and an excessive number of adaptive iteration steps. 
	
	To improve the efficiency of the AFEM, the purpose of this paper is to
	propose an efficient AFEM based on high-accuracy techniques, including the superconvergence recovery and the CVDT mesh optimization. In the adaptive procedure, the mesh optimization is embedded to preserve that the adaptive mesh is of high quality which assures the superconvergence property of the recovered gradients and the asymptotical exactness of the gradient recovery based on a posteriori error estimator. Furthermore, a tailored adaptive strategy ensures the efficiency of the whole adaptive algorithm.
	
	The remainder of the paper is organized as follows: In Section \ref{sec:model}, we provide the linear elliptic problem and its finite element scheme. Some poor performance of the standard AFEM is also presented, including the non-asymptotic accuracy of the residual-type a posteriori error estimation and the excessive iteration steps resulting from the bisection mesh refinement. We then propose the high-accuracy techniques based on AFEM in Section \ref{sec:AFEM}. Several numerical results are presented in Section \ref{sec:num} to show the effectiveness of the proposed algorithm.

	\section{Model problem and adaptive finite element method.}\label{sec:model}
	
	\subsection{Model problem and its finite element discretization.}
	Let $\Omega\subset\mathbb{R}^2$ be a polygonal domain with boundary $\partial\Omega$. Consider the following elliptic PDE with homogeneous boundary condition
	\begin{equation}\label{model}
		\left\{\begin{aligned}
			-\nabla\cdot(A\nabla u) &=f\quad \text{in}~\Omega, \\
			u&=0\quad \text{on}~\partial \Omega,
		\end{aligned}\right.
	\end{equation}
	where $f\in L^2(\Omega)$, and $A=(a_{ij})_{2\times 2}$ is a given positive definite symmetric matrix function in $L^{\infty}(\Omega)$. 
	
	Let $V$ be the Hilbert space $H_0^1(\Omega)$, defined as
	$$H_0^1(\Omega):=\{v\in H^1(\Omega):v=0  ~\text{on}~  \partial\Omega \}.$$
	The variational formulation of problem \eqref{model} reads: find $u\in H_0^1(\Omega)$ such that
	\begin{equation}\label{yuan}
		a(u,v)=(f,v) \quad\forall\; v\in H_0^1(\Omega),
	\end{equation}
	where 
	\[a(u,v)=\int_\Omega A \nabla u\cdot\nabla v{\rm d}x,
	\qquad (f,v)=\int_\Omega fv {\rm d}x.\] 
	For any $v\in V$, we define its energy norm $\|v\|_V$ associated with the bilinear form
	\[\||v|\|:=(a(v,v))^{1/2}.\]
	It is easy to verify that the bilinear form is bounded and elliptic on $V\times V$ 
	\begin{equation}\label{vfe}
		\begin{aligned}|a(u,v)|\leq &C\||u|\|\||v|\|, \qquad \forall u,v\in V\\
			a(v,v)\geq &\alpha \||v|\|^2, \qquad \qquad\forall v\in V.
		\end{aligned}    
	\end{equation}
	
	Assume that $\mathcal{T}_h$ is a conforming and shape-regular triangulation of $\Omega$. Let $h_T$ be the diameter of the element $T\in\mathcal{T}_h$, and $V_h\subset V$ be the continuous piecewise linear finite element space associated with $\mathcal{T}_h$
	\[V_h=\{v\in H_0^1(\Omega): v\in P_1(T), \forall T\in\mathcal{T}_h\},\]
	where $P_1(T)$ denotes the set of linear polynomials on $T$.
	Set $V_h^0=V_h\bigcap H_0^1(\Omega)$,
	then the finite element scheme for \eqref{yuan} is to find $u_h\in V_h^0 $ such that
	\begin{equation}\label{lisan}
		a(u_h,v_h)=(f,v_h),\qquad \forall\,v_h\,\in \,V_h^0.
	\end{equation}
	The existence and uniqueness of the finite element solution are provided by the Lax-Milgram theorem.

	\subsection{Poor performance of the standard adaptive finite element method.}
	In some cases, the solution of the problem \eqref{model} could have sharp gradients or singularity when the domain is concave or the coefficient is discontinuous. 
	AFEM has been widely used to improve the accuracy of numerical approximations to \eqref{model}. The general idea of AFEM is to refine or adjust the finite element space such that the errors are equidistributed which leads to optimal order of convergence. The standard $h$-adaptive finite element method  
	consists of loops of the form:
	\[\textbf{SOLVE}\rightarrow \textbf{ESTIMATE} \rightarrow \textbf{MARK} \rightarrow \textbf{REFINE}.\]
	The procedure \textbf{SOLVE} solves finite element equation \eqref{lisan} to obtain the numerical approximation $u_h$. 
	The procedure \textbf{ESTIMATE} calculate the element-wise error indicators $\{\eta_{h, T}\}_{T\in \mathcal{T}_h}$ and the global error estimator $\eta_h=\left(\sum\limits_{T\in \mathcal{T}_h}\eta_{h, T}^2\right)^{1/2}$. The a posteriori error estimators are an essential part of the \textbf{ESTIMATE} step. By using the
	information from the approximated solution and the known data, the
	a posteriori error estimator provides information about the size
	and the distribution of the error of the finite element
	approximation. 
	The procedure \textbf{MARK} identifies a subset $\mathcal{M}$ of elements to be refined according to the maximum marking strategy 
	or the D\"{o}rfler marking strategy.
	After choosing a set of marked elements, the procedure \textbf{REFINE} partitions the marked triangles such that the mesh obtained is still conforming and shape regular. Two well-known classes of mesh refinement algorithms are the bisection method and the red-green method.
	
	\normalem
	\begin{algorithm}
		\caption{Standard adaptive finite element algorithm. }\label{algAFEM}
		\KwIn{Domain $\Omega$, right hand side function $f$, coefficient matrix $A$, tolerance $TOL<1$, a parameter $0<\theta<1$}
		\KwOut{Sequence of mesh $\{\mathcal{T}_h^{(k)}\}$ and finite element approximations $\{u_h^{(k)}\}$}
		Set $k=0$, $\eta_h^{(k)}=1$\;
		Generate initial mesh $\mathcal{T}_h^{(k)}$ \;
		\While {$\eta_{h}^{(k)} \geq TOL$}{
			Set $k=k+1$ \;
			
			Solve the equation \eqref{lisan} on mesh $\mathcal{T}_h^{(k)}$ to get solution $u_h^{(k)}$ \;
			
			Calculate local error indicator $\{\eta_{h,T}^{(k)}\}_{T\in\mathcal{T}_h^{(k)}}$ and global error estimator $\eta_{h}^{(k)}$ \;
			
			Mark a set $\mathcal{M}_k\subset \mathcal{T}_k$ with minimum number such that \[\sum\limits_{T\in \mathcal{M}_k}\eta_{h, T}^2\geq \theta \sum\limits_{T\in \mathcal{T}_k}\eta_{h, T}^2\]
			
			Refine the elements $T\in \mathcal{M}_k$ and necessary elements by the newest vertex bisection method to a new mesh $\mathcal{T}_k$ \;
		}
	\end{algorithm}

	We formulate the adaptive algorithm with the D\"{o}rfler marking strategy and the newest vertex bisection method in Algorithm \ref{algAFEM}.
	Starting from an initial mesh $\mathcal{T}_0$, Algorithm \ref{algAFEM} produces a sequence of meshes $\{\mathcal{T}_k\}_{k\in\mathbb{N}_0}$ and finite element solutions $\{u_k\}_{k\in\mathbb{N}_0}$. Under some modest assumptions, we can derive the contraction of error estimator and the total error between two consecutive levels $\mathcal{T}_k$ and $\mathcal{T}_{k+1}$, then can prove AFEM will terminate in a finite number of steps for a given tolerance $tol$ and yield a convergent approximation $u_k$ on an adaptive mesh $\mathcal{T}_k$. 
	
	In the adaptive finite element methods, a posteriori error estimators are able to locate accurately sources of global and local error in the approximation. 
	Mainly, there are three types of a posteriori error estimator include: residual type \cite{br, Babu_1978, Bernardi2000}, recovery type \cite{cz, yan2001, fa, zz}, solving auxiliary problem \cite{AD2000}.

	Denote the set of all edges of the triangulation $\mathcal{T}_h$ as
	$$\mathcal{E}:=\mathcal{E}_\Omega\cup\mathcal{E}_{\partial\Omega},$$
	where $\mathcal{E}_\Omega$ is the set of all interior element edges and $\mathcal{E}_{\partial\Omega}$ is the sets of boundary edges. For each $e\in \mathcal{E}_{\Omega}$, we fixed a unit norm vector $\textbf{n}_e$. Let $T_e^+$ and $T_e^-$ be two triangles sharing the edge $e$. 
	For the model problem \eqref{model} and its finite element scheme \eqref{lisan}, the standard residual type estimator \cite{MKN2005} on element $T$ is defined as
	\begin{equation}\label{residual_T}
		\eta_{T,res}^2=h_T^2\|R_T(u_h)\|_{0,T}^2+	\sum_{e\in \partial T} h_e\|J_e(A \nabla u_h)\|_{0,e}^2,
	\end{equation}
	where the element residual $R_T(u_h)$ and the jump of flux are defined as 
	\[\begin{aligned}
		R_T(u_h):=f+\nabla\cdot(A\nabla u_h)\qquad & \text{in}\ T\in\mathcal{T}_h,\\
		J_e(A \nabla u_h):=[A \nabla u_h\cdot \textbf{n}_e]=(A \nabla u_h|_{T_e^+}-A \nabla u_h|_{T_e^-})\cdot \textbf{n}_e\qquad &  \text{on}\ e\in\mathcal{E}_\Omega,\end{aligned}\]
	and $J_e(A \nabla u_h)=0$	for any $e\in\mathcal{E}_{\partial\Omega}$. Correspondingly, the global residual type estimator is defined by
	\begin{equation}\label{residual}
		\eta_{res}^2=\sum_{T\in \mathcal{T}_h}h_T^2\|R_T(u_h)\|_{0,T}^2+	\sum_{e\in \mathcal{E}} h_e\|J_e(A \nabla u_h)\|_{0,e}^2.
	\end{equation}
	
	\begin{thm}[\cite{MKN2005}]\label{res_1}
		There exists constant $C$ depending on the domain $\Omega$, the coefficient function $A$, and the regularity of $\mathcal{T}_h$ such that the residual type estimator is the globally upper and locally lower bounds of the error's energy norm, i.e.,
		\[||| u-u_h|||\leq C\eta_{res}, \quad \eta_{T,res}^2\le C\big(|||u-u_h|||_{\omega_T}^2 + osc_h^2(f,\omega_T)\big),\]
		where $osc(f)$ is the high-order oscillation term, which is defined by
		\[osc_h^2(f,\omega_T)=\sum_{T\in\omega_T}h_T^2\|R_T(u_h)-\overline{R_T(u_h)}\|_{0,T}^2.\] 
		$\overline{R_T(u_h)}$  denotes the element average of $R_T(u_h)$ on $T$, and $\omega_T$ represents the set of elements that share at least one side with $T$.
	\end{thm} 
	
	Theorem \ref{res_1} shows the reliability and efficiency of the residual estimator. 
	Although the adaptive algorithm derived by this type of estimator is quasi-optimal in the sense that the adaptive meshes generated by this method provide the highest possible convergence rate, it is not
	favorable in some aspects. For one, the bounds provided by the error estimator contain unknown constants, which depend on the mesh quality,  the polynomial degree used in the finite element space, and the known data of the considered problem. 
	
	In the following numerical test, we shall show that the adaptive finite element algorithm driven by this estimator may not stop timely, and thus waste a large amount of computational cost.
	
	\begin{example}\label{eg1}
		Consider the problem with geometric singularity
		\begin{equation}
			\left\{\begin{aligned}\label{egL}
				-\Delta u=0&\quad \text{in}~\Omega, \\
				u=g&\quad \text{on}~\partial \Omega,
			\end{aligned}\right.
		\end{equation}
		where the L-shape domain $\Omega= (-1,1)^2\backslash (0,1)\times(-1,0)$. In this example, we choose the exact solution $u$ as
		\[u(x,y)=r^{2/3}\sin(2\theta/3),\quad r=\sqrt{x^{2}+y^{2}},\quad \theta=\tan^{-1}(y/x),\]
		and $g$ is determined from $u$. The numerical solution has a singularity at the origin point $(0,0)$.
	\end{example}
	
	\begin{figure}[ht]
		\begin{minipage}{0.29\linewidth}
			\centerline{\includegraphics[width=5.5cm]{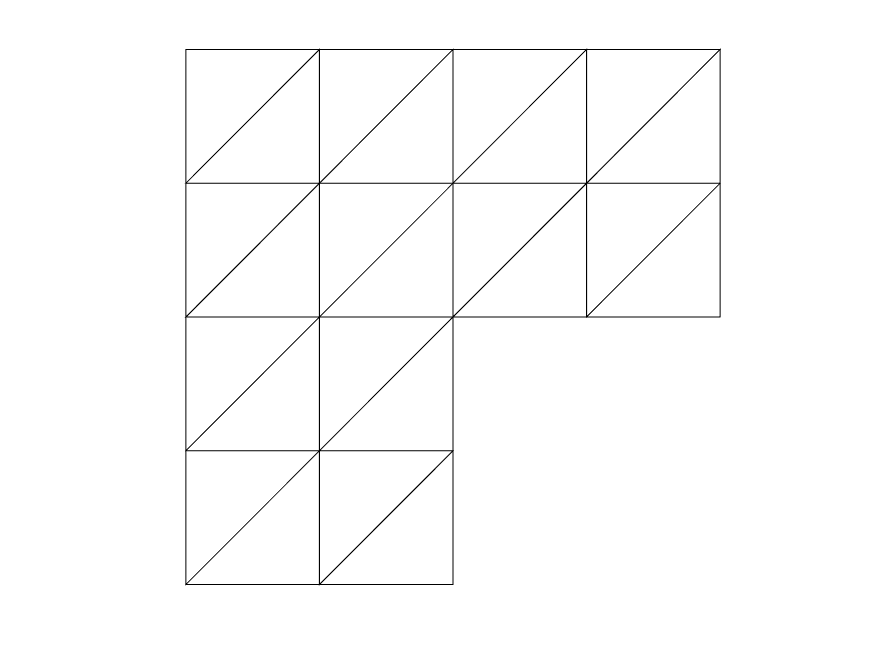}}
			\centerline{(a)}
		\end{minipage}
		\begin{minipage}{0.31\linewidth}
			\centerline{\includegraphics[width=5.5cm]{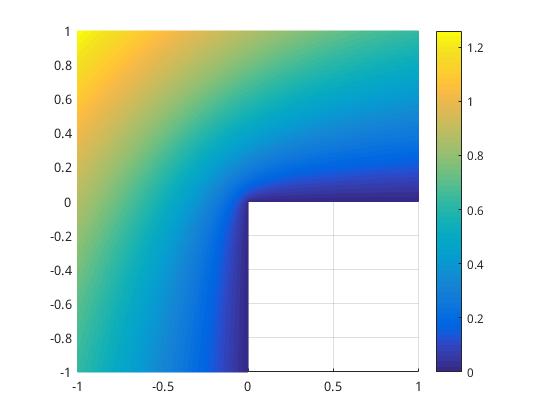}}
			\centerline{(b)}
		\end{minipage}
		\begin{minipage}{0.3\linewidth}
			\centerline{\includegraphics[width=5.5cm]{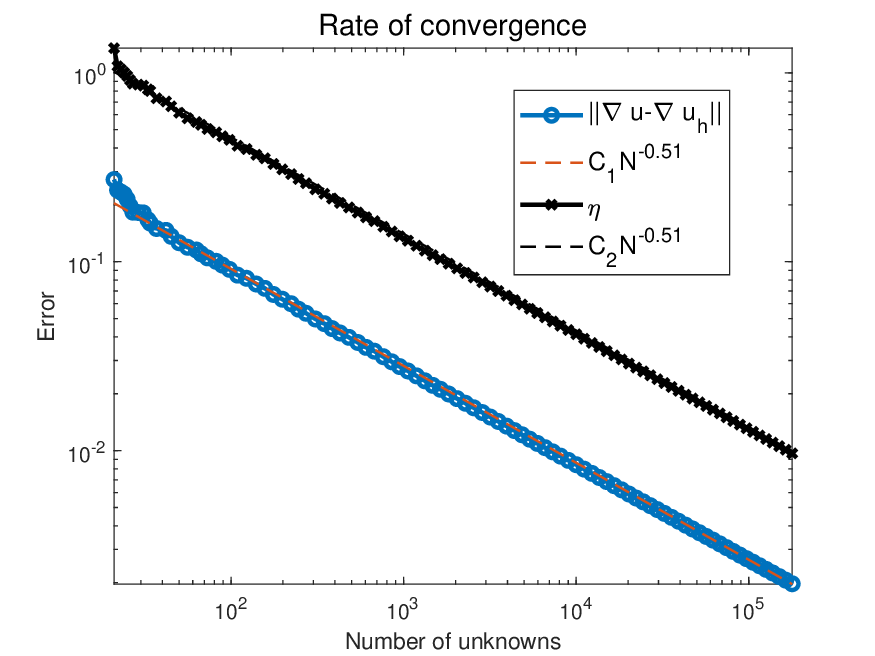}}
			\centerline{(c)}
		\end{minipage}
		\caption{Example \ref{eg1}, (a) initial mesh; (b) numerical solution; (c)
			errors.}\label{figR1}
	\end{figure}
	
	\begin{figure}[!ht]
		\begin{minipage}{0.3\linewidth}
			\centerline{	\includegraphics[width=5.5cm]{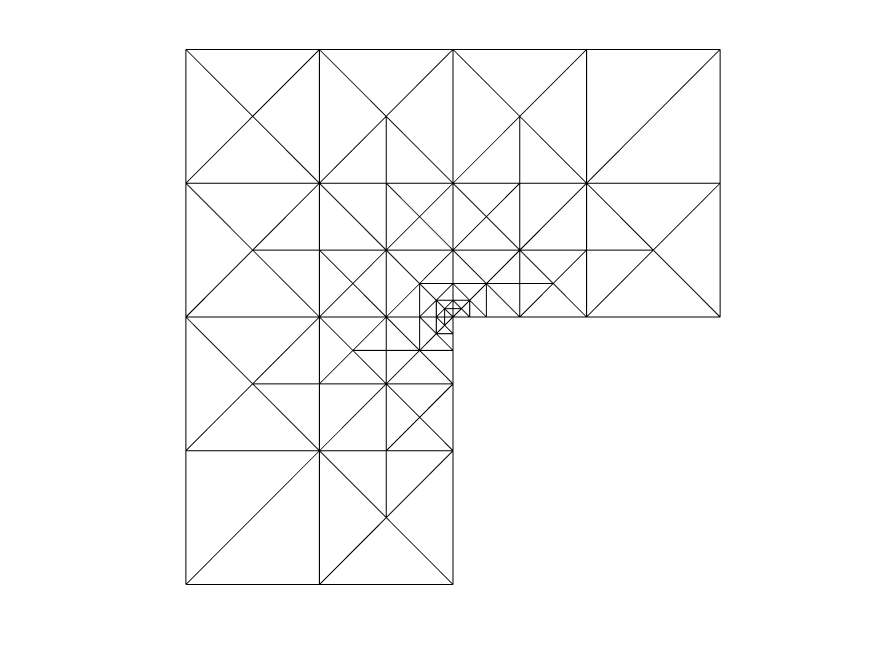}}
		\end{minipage}
		\begin{minipage}{0.3\linewidth}
			\centerline{	\includegraphics[width=5.5cm]{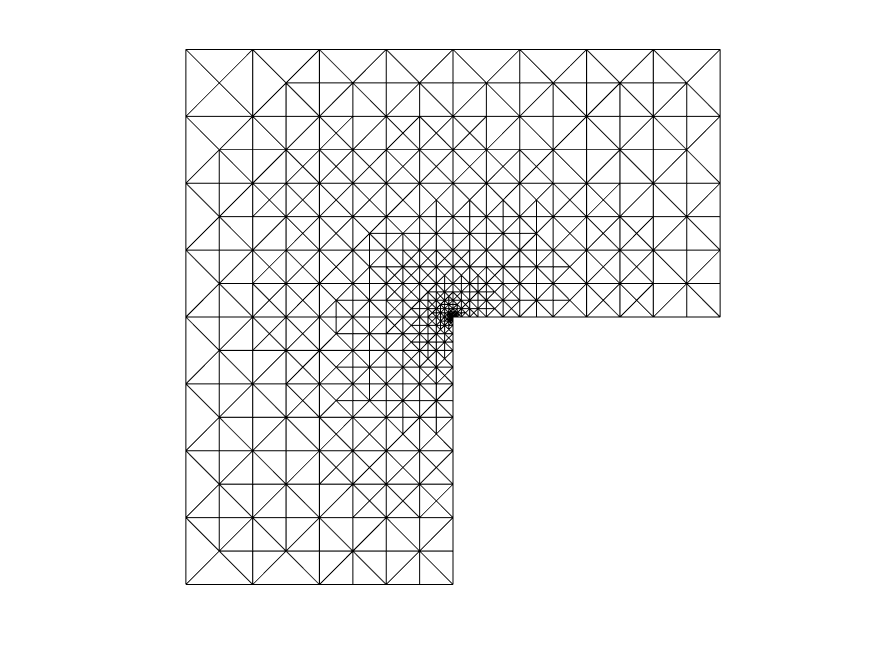}}
		\end{minipage}\begin{minipage}{0.3\linewidth}
			\centerline{	\includegraphics[width=5.5cm]{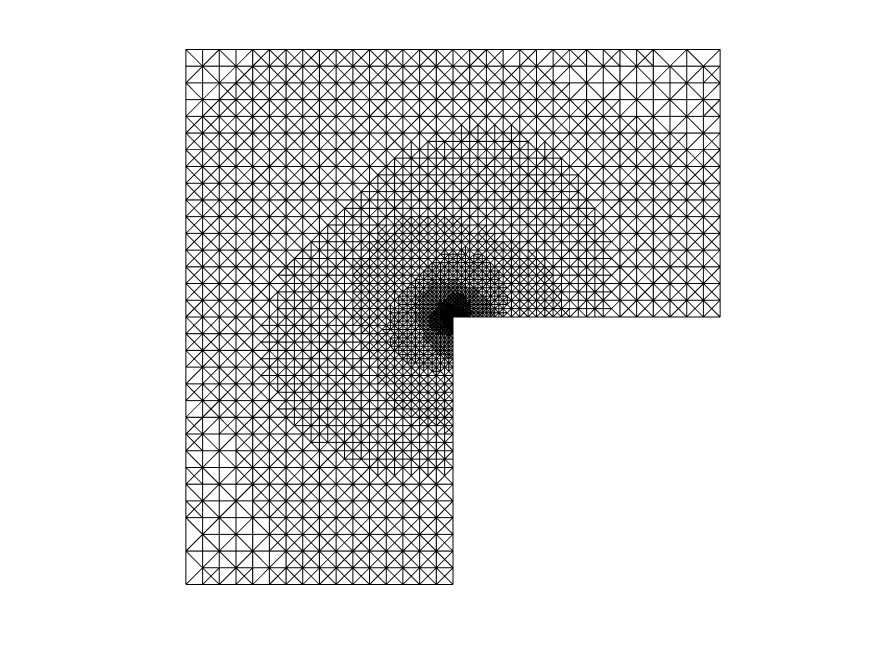}}
		\end{minipage}
		
		\begin{minipage}{0.3\linewidth}
			\centerline{	\includegraphics[width=5.5cm]{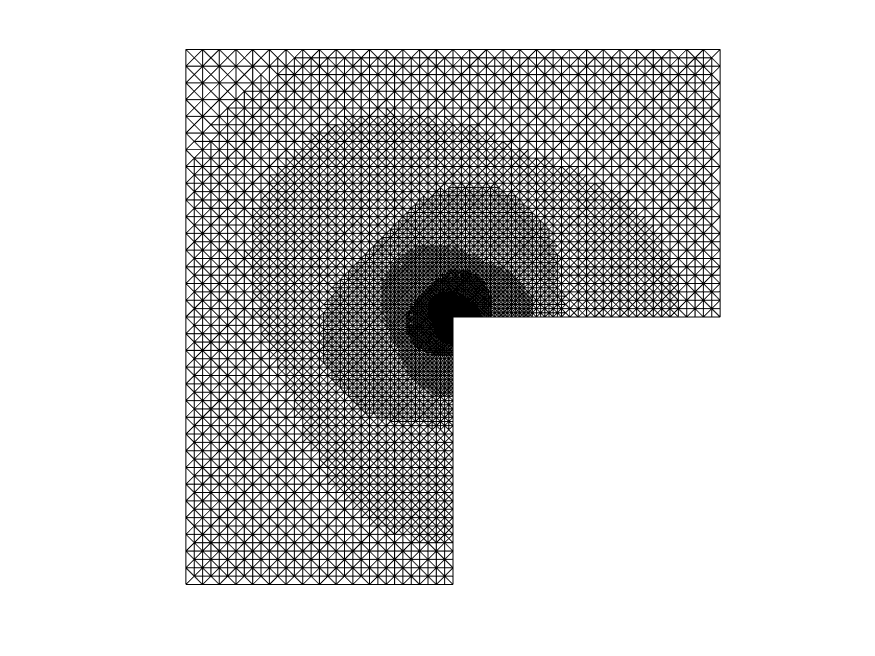}}
		\end{minipage}
		\begin{minipage}{0.3\linewidth}
			\centerline{	\includegraphics[width=5.5cm]{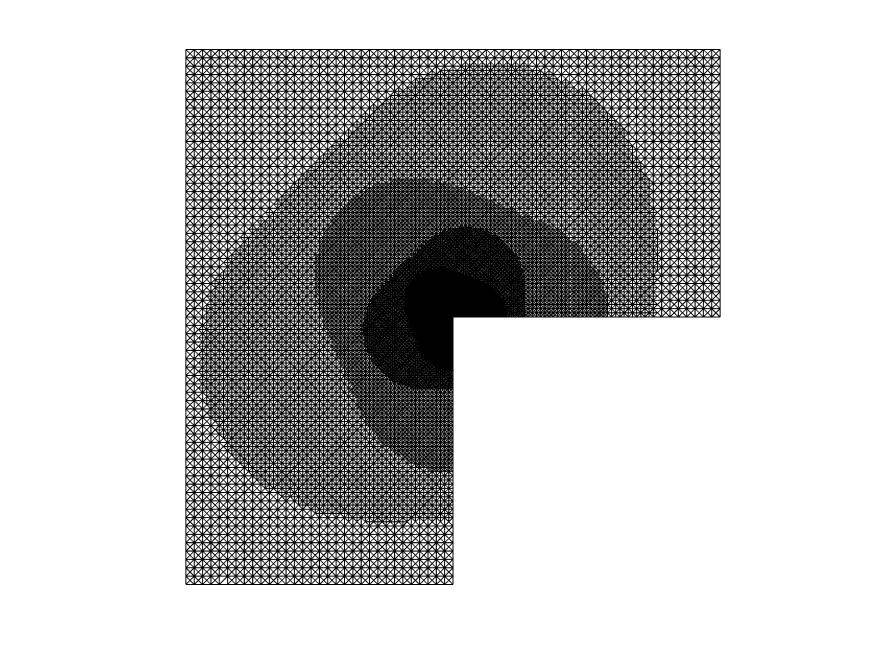}}
		\end{minipage}
		\begin{minipage}{0.3\linewidth}
			\centerline{	\includegraphics[width=5.5cm]{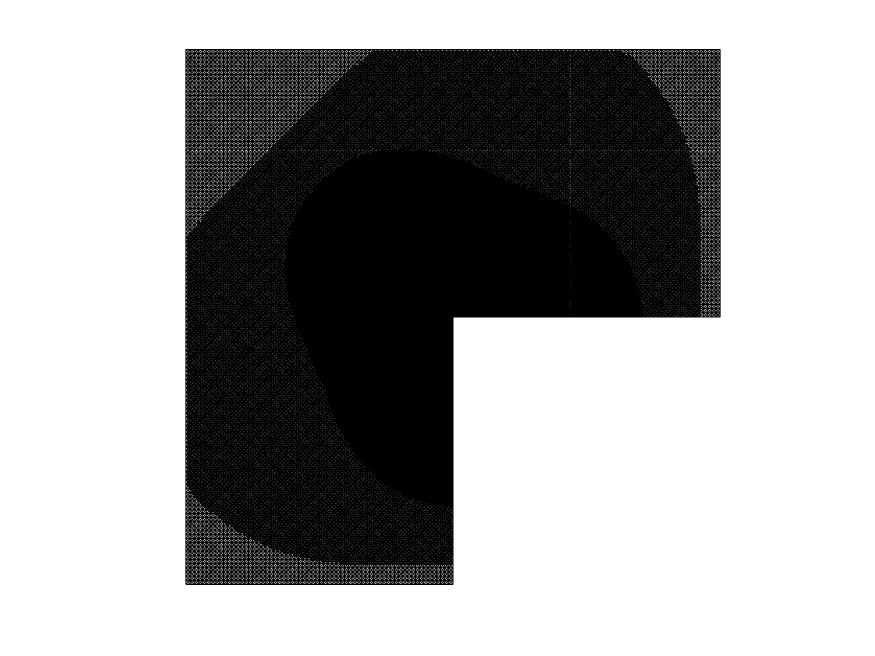}}
		\end{minipage}
		\caption{Example \ref{eg1}, adaptive meshes at steps: 10, 20, 30, 35, 40 and 50.}\label{figR2}
	\end{figure}
	
	\begin{table}[ht]
		\caption{Example \ref{eg1}, data of mesh vertices, error, and error estimator.}\label{Lshpaedata} \centering
		\begin{tabular}{cccc} \hline
			$k$    &    $N$       & $||\nabla u-\nabla u_h||$ & $\eta_{res}$                  \\ \hline
			1     &     21       & $2.7205e-01$  & $1.3499e+00$ \\
			$\vdots$   &     $\vdots$      & $\vdots$                       & $\vdots$                      \\
			34    &     6781     & $1.0225e-02$  & $5.0067e-02$ \\
			\bf{35}    & \bf{8073}   & $\bf{9.3704e-03}$  & $\bf{4.5878e-02}$ \\
			36    &     9661     & $8.5595e-03$  & $4.1932e-02$ \\
			$\vdots$   &     $\vdots$      & $\vdots$                       & $\vdots$                      \\
			52    &     135971   & $2.2638e-03$  & $1.1088e-02$ \\
			\bf{53} &\bf{159679} & $\bf{2.0876e-03}$ & $\bf{1.0227e-02}$ \\\hline
		\end{tabular}
	\end{table}
	
	We apply Algorithm \ref{algAFEM} to solve \eqref{egL} with $TOL=0.01$ and $\theta=0.3$. The initial mesh, numerical solution, and the convergence history of error and estimator are plotted in Figure \ref{figR1}. The adaptive meshes are presented in Figure \ref{figR2}. 
	The numerical results show that the mesh refinement process meets the singularity of the solution at the reentrant corner, and the convergence rates of the errors is quasi-optimal. 
	
	However, the residual type estimator is approximately $5$ times the exact error, which in turn does not terminate the adaptive procedure when the gradient error $\|\nabla u-\nabla u_h\|<TOL$. To see this clearly, we list the data of mesh vertices, error $\|\nabla u-\nabla u_h\|$ and the residual type estimator $\eta_{res}$ in Table \ref{Lshpaedata}. We see clearly that the adaptive iteration terminates until the adaptive mesh arrives at $159679$ vertices, and the corresponding gradient error is $2.0876e-03$. It is important to note that the error satisfies $\|\nabla u-\nabla u_h\|<TOL$ at the $35$th step with the adaptive mesh has $8073$ vertices, which means that the adaptive iterations after the $35$th iteration are unnecessary. Furthermore, these redundant iteration steps consume much more computational cost. The unsatisfactory over-refinement is due to the non-asymptotically exact residual type estimator. In addition, beginning with the initial mesh with $21$ vertices, the mesh adaption process $34$ steps to obtain a mesh with $8073$ vertices. The reason for this phenomenon is that the bisection method is a layer-by-layer refinement method, which makes it difficult to achieve the target mesh size in just a few steps. 
	
	Due to the use of non-asymptotically exact a posteriori error estimate and ineffective mesh adaptive strategy, the standard AFEM usually requires many adaptive iteration steps, and thus reduces the computational efficiency of the adaptive algorithms.  
	Fortunately, we can use the recovery type estimator to eliminate one of the shortcomings of the standard residual type estimator.  
	The gradient recovery-based a posteriori error estimator, which uses a certain norm of the difference between the direct and post-processed approximations of the gradient as an indicator, is asymptotically exact when the recovered gradient is superconvergent.
	Gradient recovery is a post-processing technique
	that reconstructs improved gradient
	approximations from finite element solutions as well as to explore their use in adaptive computations \cite{cz, yan2001, fa,  hy, hy1, LXYC2024, yz, za}.
	Denote the set of all mesh nodes of $\mathcal{T}_h$ by $\mathcal{N}_h$. 
	For the finite element solution $u_h\in V_h$, denote the recovered gradient space $W_h=V_h\times V_h$ and the gradient recovery operator $G: V_h\rightarrow W_h$. 
	For each node $z\in\mathcal{N}_h$, we first define the recovered gradient $G(\nabla u_h(z))$ by the least square fitting or projection methods on the element patch $\omega_z$,
	then obtain the recovered gradient $G(\nabla u_h)$ on the whole domain by interpolation
	\[G(\nabla u_h)=\sum_{z\in\mathcal{N}}G(\nabla u_h(z))\phi_z,\]
	where $\phi_z$ is the Lagrange basis of finite element space $V_h$ associated with $\mathcal{T}_h$.
	
	The basic principle behind gradient recovery-based error estimator is to apply some inexpensive post-processing to the gradient of the finite element approximation, $\nabla u_h\rightarrow G(\nabla u_h)$, so
	that the recovered gradient $G(\nabla u_h)$ provides a better estimate of the true gradient $\nabla u$
	than $\nabla u_h$ does. Then define the local a posteriori error estimator as
	\begin{equation}\label{etarec}
		\eta_{T, rec}^2=\|G(\nabla u_h)-\nabla u_h\|_{0, T},\quad \eta_{rec}=\|G(\nabla u_h)-\nabla u_h\|.\end{equation}
	If the recovery procedure can improve the approximation in the sense that 
	\[\|G(\nabla u_h)-\nabla u\|\leq \beta\|\nabla u-\nabla u_h\|,\qquad {\rm with}\ 0\leq\beta<1,\]
	we immediately obtain the upper and lower bounds for the recovery-based error estimator
	\[\frac{1}{1+\beta}\|G(\nabla u_h)-\nabla u_h\|
	\leq \|\nabla u-\nabla u_h\|\leq \frac{1}{1-\beta}\|G(\nabla u_h)-\nabla u_h\|.
	\]
	In addition, if the recovery method is superconvergent, 
	\[\|G(\nabla u_h)-\nabla u\|=o(\|\nabla u-\nabla u_h\|),\]
	this means $\beta\rightarrow 0$ and  
	\[\frac{\|G(\nabla u_h)-\nabla u_h\|}{\|\nabla u-\nabla u_h\|}=1+o(1),\]
	then the corresponding error estimator is asymptotic exactness.
	
	The following theorem gives the reliability and efficiency of the gradient recovery type error estimator.
	\begin{thm}[\cite{yan2001,LXYC2024}]
		There exist constants  $C_1$ and $C_2$  such that the gradient recovery type estimator satisfies
		$$C_1||| u-u_h||| - \varepsilon_1\leq \eta_1\leq C_2||| u-u_h|||+\varepsilon_2,$$
		where $\varepsilon_1$ and $\varepsilon_2$ are high-order terms.
	\end{thm}
	
	\begin{example}\label{eg3-1}
		We consider the problem in Example \ref{eg1}, and use the recovery type error estimator \eqref{etarec} and the newest vertex bisection refinement in Algorithm \ref{algAFEM} to compare the numerical performance of the residual estimator in Example \ref{eg1}. 
	\end{example}
	
	\begin{figure}[!ht]
		\begin{minipage}{0.3\linewidth}
			\centerline{	\includegraphics[width=5.5cm]{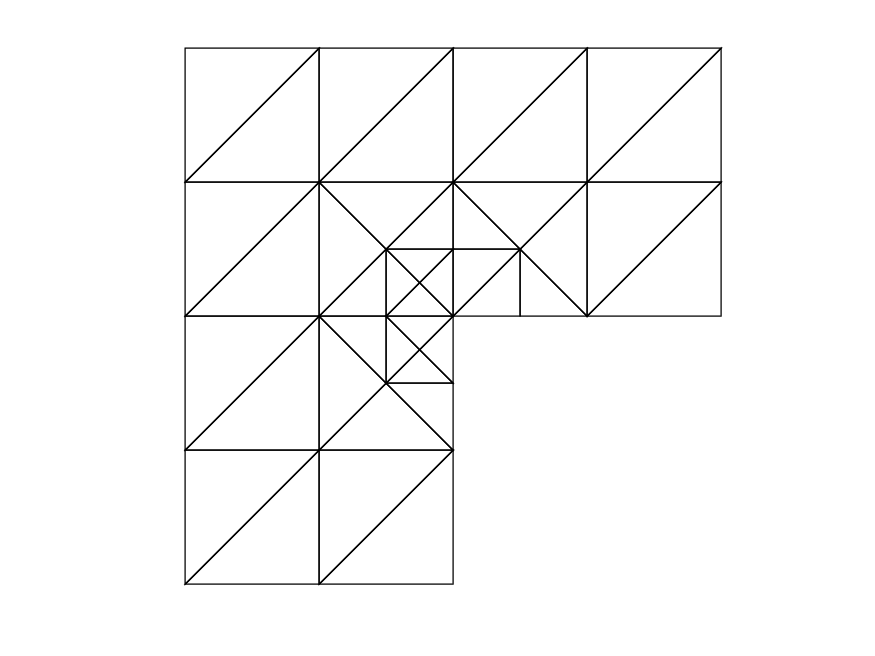}}
		\end{minipage}
		\begin{minipage}{0.3\linewidth}
			\centerline{	\includegraphics[width=5.5cm]{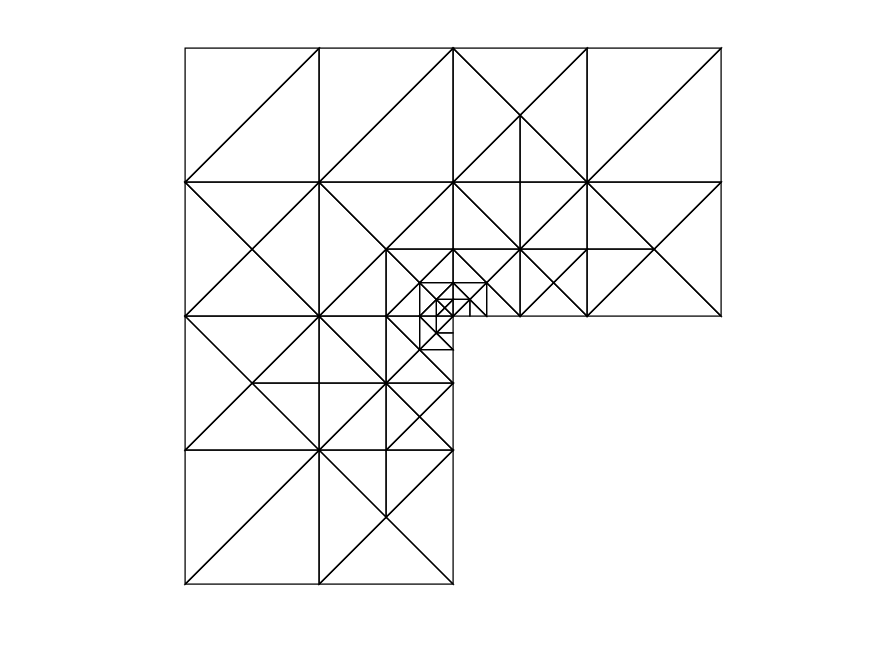}}
		\end{minipage}
		\begin{minipage}{0.3\linewidth}
			\centerline{	\includegraphics[width=5.5cm]{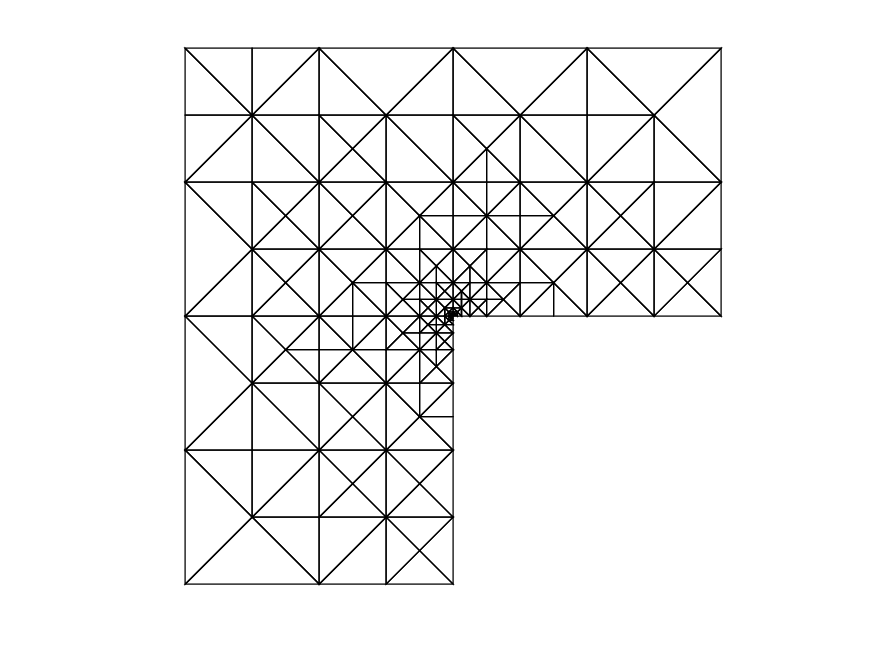}}
		\end{minipage}
		
		\begin{minipage}{0.3\linewidth}
			\centerline{	\includegraphics[width=5.5cm]{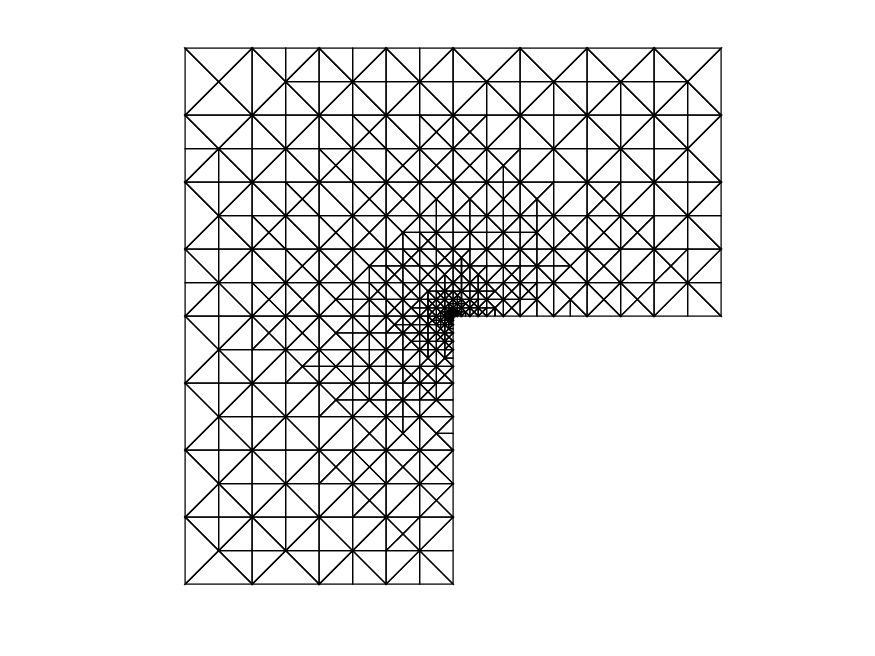}}
		\end{minipage}
		\begin{minipage}{0.3\linewidth}
			\centerline{	\includegraphics[width=5.5cm]{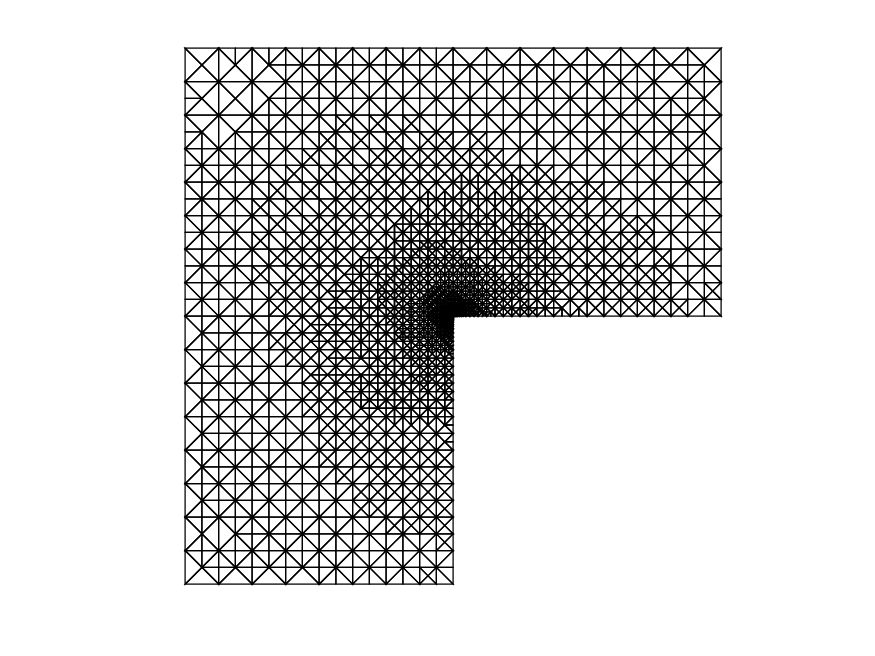}}
		\end{minipage}
		\begin{minipage}{0.3\linewidth}
			\centerline{	\includegraphics[width=5.5cm]{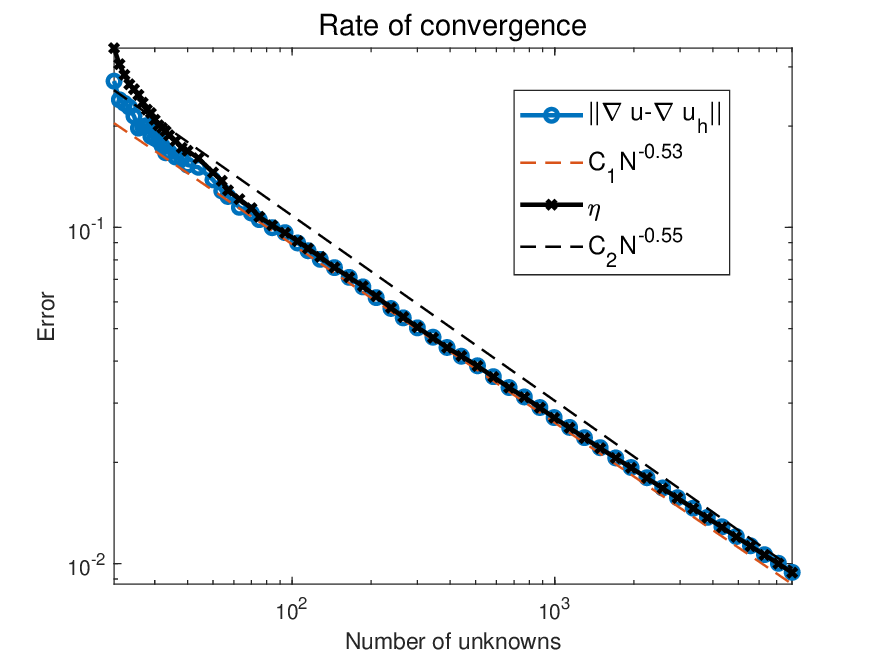}}
		\end{minipage}
		\caption{Example \ref{eg3-1}, 
			adaptive meshes of steps 10, 20, 30, 40, 50, and errors.}\label{figLSCT}
	\end{figure}
	
	\begin{table}[ht]
		\caption{Example \ref{eg3-1}, data of mesh vertices, error and the recovery type error estimator.}\label{recordata} \centering
		\begin{tabular}{cccc} \hline
			$k$    &    $N$       & $||\nabla u-\nabla u_h||$ & $\eta_{rec}$                  \\ \hline
			1     &     21       & $2.7205e-01$  & $1.3499e+00$ \\
			$\vdots$   &     $\vdots$      & $\vdots$                       & $\vdots$                      \\
			30	&	145	&	$7.5880e-02$	&	$7.5964e-02$\\
			31	&	165	&	$7.0976e-02$	&	$7.1091e-02$ \\
			32	&	186	&	$6.6507e-02$	&	$6.6820e-02$ \\
			$\vdots$   &     $\vdots$   & $\vdots$                       & $\vdots$                      \\
			59	&	6289&	$1.0633e-02$	&	$1.0613e-02$ \\
			60	&	7105&   $1.0023e-02$	&	$1.0001e-02$ \\
			\bf{61} &\bf{8018} & $\bf{9.4224e-03}$ & $\bf{9.4003e-03}$ \\\hline
		\end{tabular}
	\end{table}
	
	Similar to Example \ref{eg1}, we take $TOL=0.01$ and the initial mesh is shown in Figure \ref{figR1} (a). The adaptive meshes and the convergence history of errors and recovery type estimators are plotted in Figure \ref{figLSCT}. Compared with Figure \ref{figR1} (c), Figure \ref{figLSCT} shows that the recovery type estimator is asymptotically exact. 
	From the data listed in Table \ref{recordata}, we can see that the adaptive loop terminates when the adaptive mesh reaches $8018$ vertices, while the results obtained by the residual type a posteriori error estimator in Table \ref{Lshpaedata} show that the adaptive mesh should reach $159679$ vertices to terminate the loops.  The efficiency of the AFEM improves a lot with the use of the recovery type estimator, while the iterations do not reduce mainly caused by the use of the bisection refinement method. To further improve the efficiency of the adaptive algorithm, we should develop a new mesh adaption method which can control the mesh size.

	\section{High accuracy techniques based adaptive finite element method}\label{sec:AFEM}
	Superconvergence techniques and mesh adaption are effective ways to obtain high-accuracy finite element approximations. Superconvergence of 
	finite element method is sensitive to the symmetry structure of the mesh, and the existing superconvergence results usually require the mesh to satisfy some strong conditions, such as uniform mesh, strongly regular mesh et al.
	Adaptive methods are now widely used in numerical partial differential equations to achieve better
	accuracy with a quasi-optimal degree of freedom. In mesh adaption, local bisection refinement would break the symmetry structure of the mesh. 
	These observations show that the superconvergence techniques and local mesh adaption may be mutually exclusive. We aim to develop new adaptive technologies that can simultaneously leverage the advantages of superconvergence techniques and mesh adaption. 
	Based on high-quality mesh optimization methods, we combine finite element superconvergence techniques with a tailored mesh adaption strategy to propose high-accuracy techniques based adaptive finite element method. In the process of adaptive solving, mesh optimization generates high-quality mesh on which gradient recovery is superconvergent, which in turn gives the recovery based a posterior error estimator asymptotically exact, and then guides the mesh adaption efficiently. 
	Furthermore, we introduce a tailored mesh adaptive strategy, based on which a high-quality mesh with a target number of vertices can be generated, thus reducing the number of adaptive iteration steps. 
	The compatibility of superconvergence gradient recovery, mesh optimization, and local adaptive refinement ensures the efficiency of the adaptive solving process.
	
	\subsection{High-quality mesh generation technique}
	The accuracy of finite element approximation depends on the sizes and shapes of the elements \cite{HQW2008, HWY2017}. In \cite{HQW2008}, superconvergence was found for the linear finite element solution on a general two-dimensional domain due to the high quality of the centroidal Voronoi Delaunay triangulation (CVDT) mesh. 
	Several robust and efficient algorithms have been developed for high-quality mesh generation and optimization \cite{chen2004, du2002, DW2005, KOKO2015650}.  In this work, we focus on centroidal Voronoi tessellation (CVT) based mesh generation and optimization algorithm \cite{Du1999, du2002}, which can generate an unstructured mesh with many desirable features, such as errors are equidistributed over the elements.
	The generating points of CVT are the mass centroids of the corresponding Voronoi regions with respect to a given density function. The localized Lloyd iteration method \cite{DW2005}, which constructs CVT by iteratively moving generators to the mass centers of Voronoi regions, effectively reduces the global distortion of element shape and sizing. The corresponding CVDT provides a high-quality unstructured mesh, on which the superconvergence property of the finite element approximation is assured. 
	If the superconvergence property of the recovered gradient holds for meshes generated in the adaptive procedure, then the gradient recovery based a posteriori error estimator is asymptotically exact. In the following, we present two numerical examples to show the efficiency of the CVDT mesh generation algorithm, and the high accuracy of finite element approximations on both uniform mesh and local adaptive refined mesh based on CVDT.
	
\begin{figure}[!ht]
		\begin{minipage}{0.45\linewidth}
			\centerline{	\includegraphics[width=6cm]{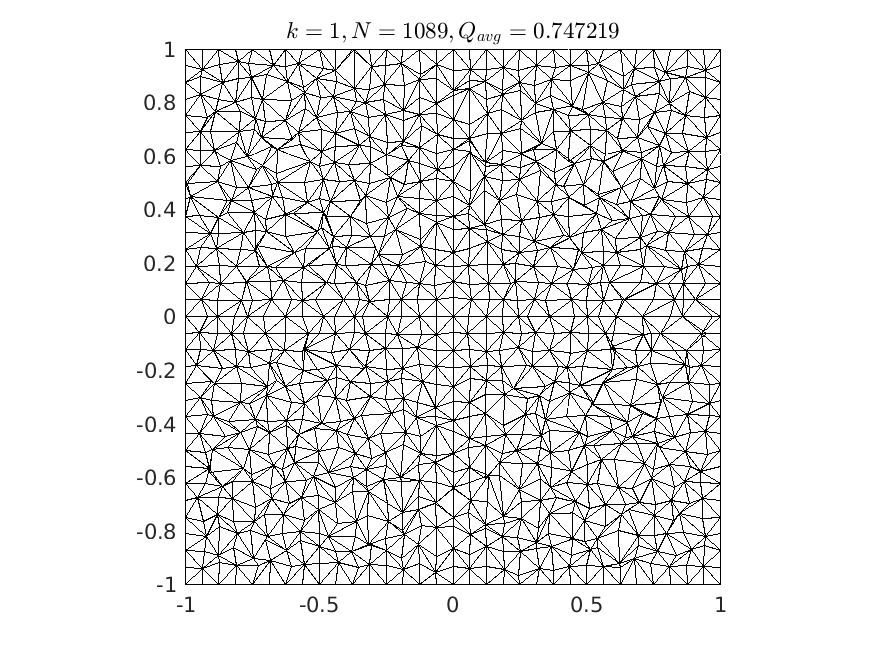}}
		\end{minipage}
		\begin{minipage}{0.45\linewidth}
			\centerline{	\includegraphics[width=6cm]{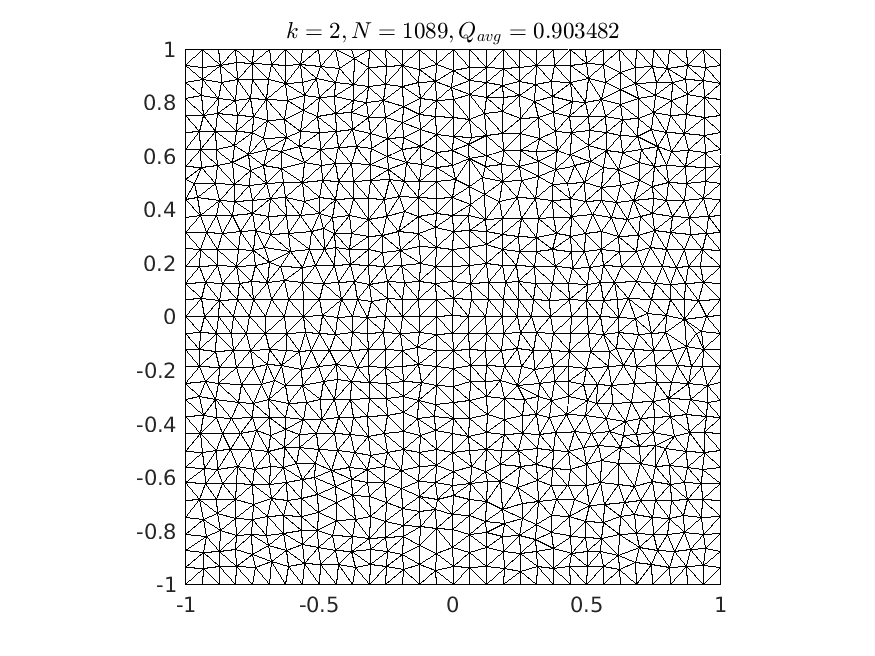}}
		\end{minipage}
		
		\begin{minipage}{0.45\linewidth}
			\centerline{	\includegraphics[width=6cm]{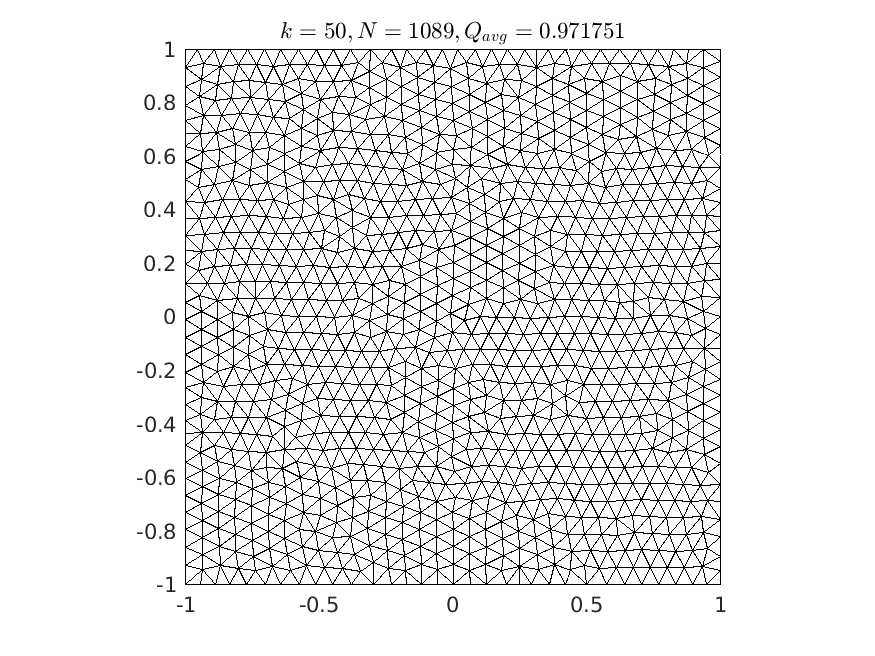}}
		\end{minipage}
		\begin{minipage}{0.45\linewidth}
			\centerline{	\includegraphics[width=6cm]{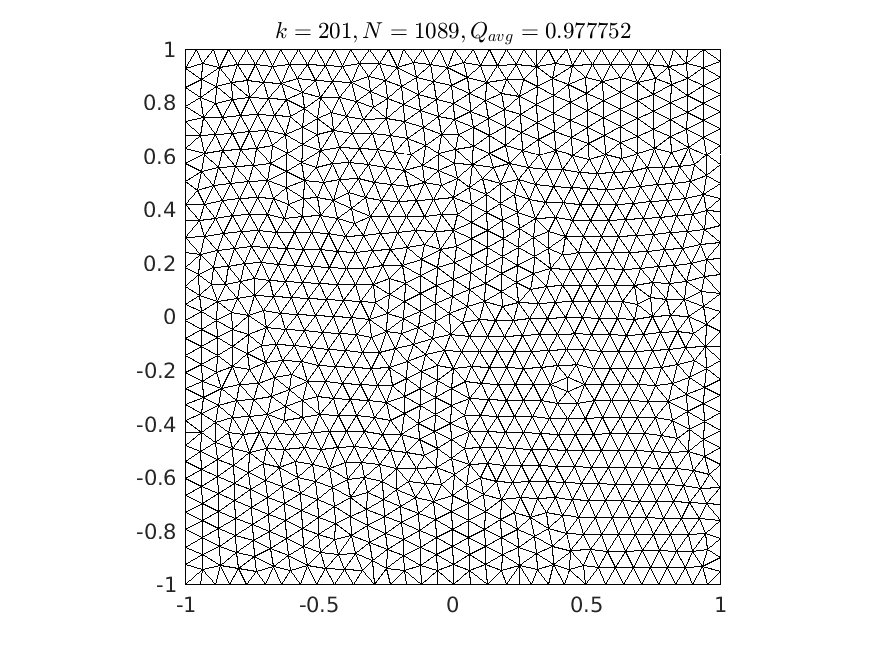}}
		\end{minipage}
		\caption{Example \ref{uniform-cvdt}, mesh generated in optimization steps $1$, $2$, $50$ and $201$.}\label{figLSCT1}
	\end{figure}

\begin{figure}[!ht]
		\centerline{	\includegraphics[width=6cm]{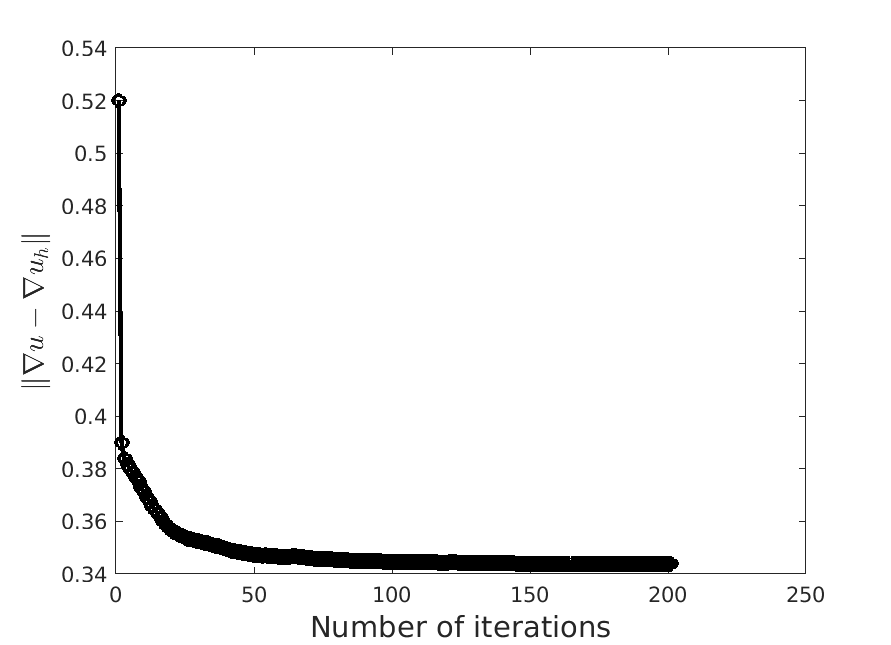}}
		\caption{Example \ref{uniform-cvdt}, history of error $\|\nabla u-\nabla u_h\|$.}\label{figLSCTerror}  
	\end{figure}

	\begin{example}\label{uniform-cvdt}
		Consider the problem
		\begin{equation}
			\left\{\begin{aligned}\label{eg0}
				-\Delta u=f&\quad \text{in}~\Omega, \\
				u=g&\quad \text{on}~\partial \Omega,
			\end{aligned}\right.
		\end{equation}
		where $\Omega=[0,1]\times[0,1]$.  The exact solution $u$ is taken as
		\[u(x,y)=cos(\pi x)cos(\pi y),\]
		and $f$ and  $g$ are determined by $u$.
		
		We randomly generate $1089$ points and obtain the initial mesh by the Delaunay algorithm, then apply the CVDT mesh algorithm to improve the mesh quality. In each mesh optimization step, we numerically solve equation \eqref{eg0} by linear finite element method and calculate the gradient error $\|\nabla u-\nabla u_h\|$. The whole mesh optimization will perform $201$ steps. 
		Figure \ref{figLSCT1} plots the mesh generated at the optimization steps $1, 2, 50$ and the final step $201$, 
		which shows that the mesh quality increases in the process of mesh optimization.
		Figure \ref{figLSCTerror} reports the history of error $\|\nabla u-\nabla u_h\|$, we see quite clearly $\|\nabla u-\nabla u_h\|$ decreases as mesh quality improves, and the error decrease fast in first few iteration steps. Therefore, in the practice calculation, one can use a small number of optimization steps to generate high quality mesh for the finite element method. 
	\end{example}
	
\begin{figure}[!ht]
		\begin{minipage}{0.3\linewidth}
			\centerline{\includegraphics[width=6cm]{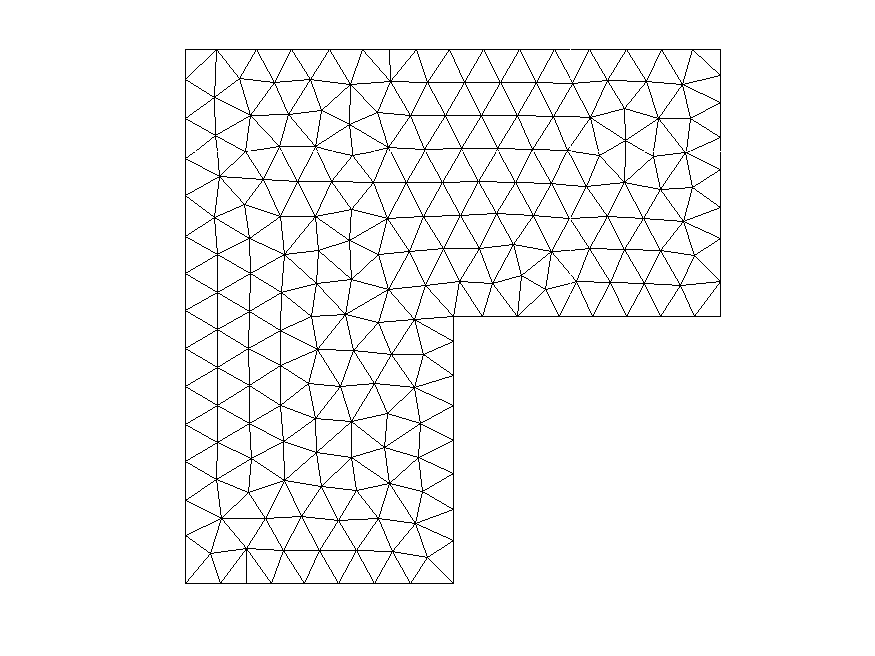}}
		\end{minipage}
		\begin{minipage}{0.3\linewidth}
			\centerline{\includegraphics[width=6cm]{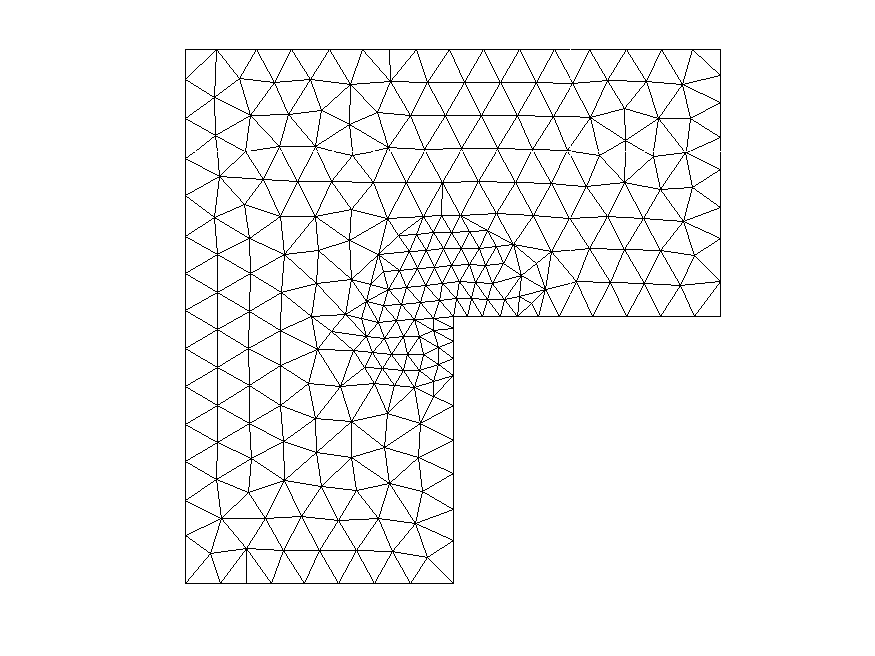}}
		\end{minipage}
		\begin{minipage}{0.3\linewidth}
			\centerline{\includegraphics[width=6cm]{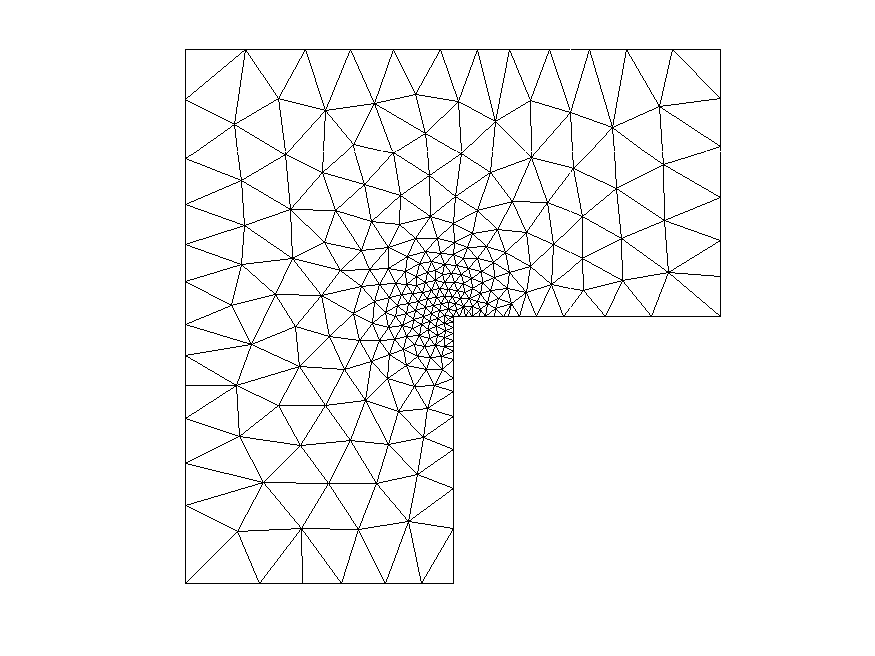}}
		\end{minipage}
		\begin{minipage}{0.3\linewidth}
			\centerline{\includegraphics[width=6cm]{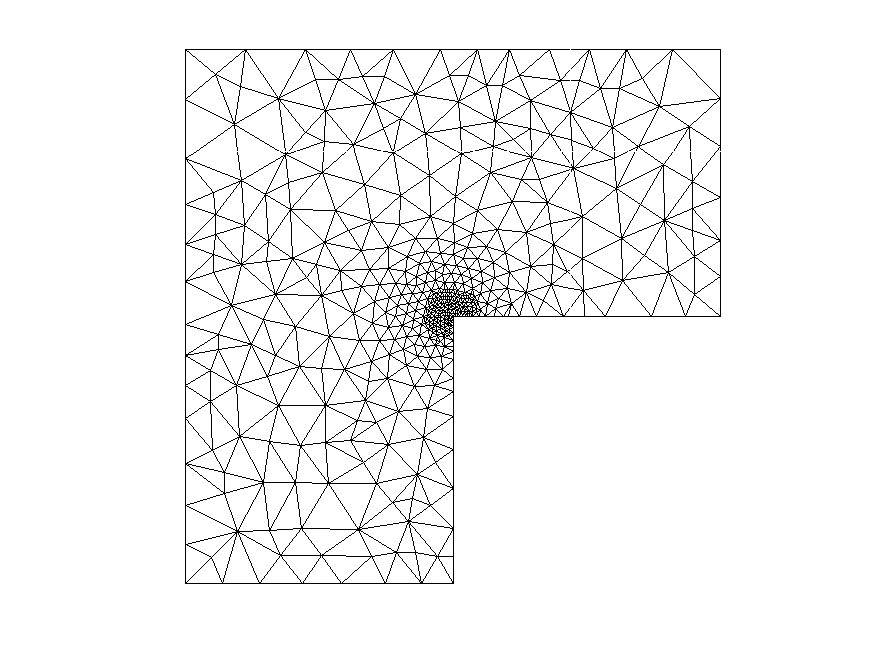}}
		\end{minipage}
		\begin{minipage}{0.3\linewidth}
			\centerline{\includegraphics[width=6cm]{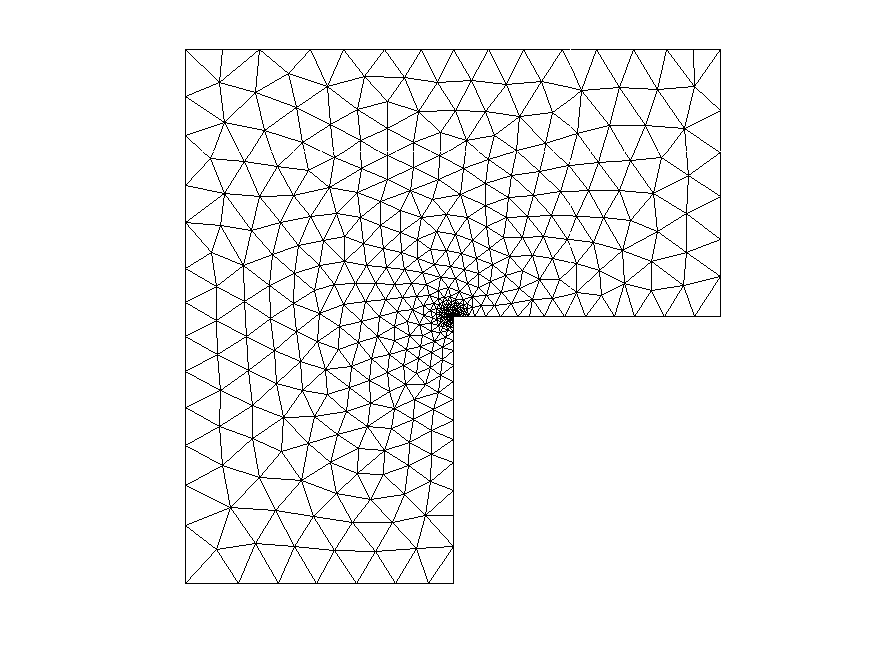}}
		\end{minipage}
		\begin{minipage}{0.3\linewidth}
			\centerline{\includegraphics[width=6cm]{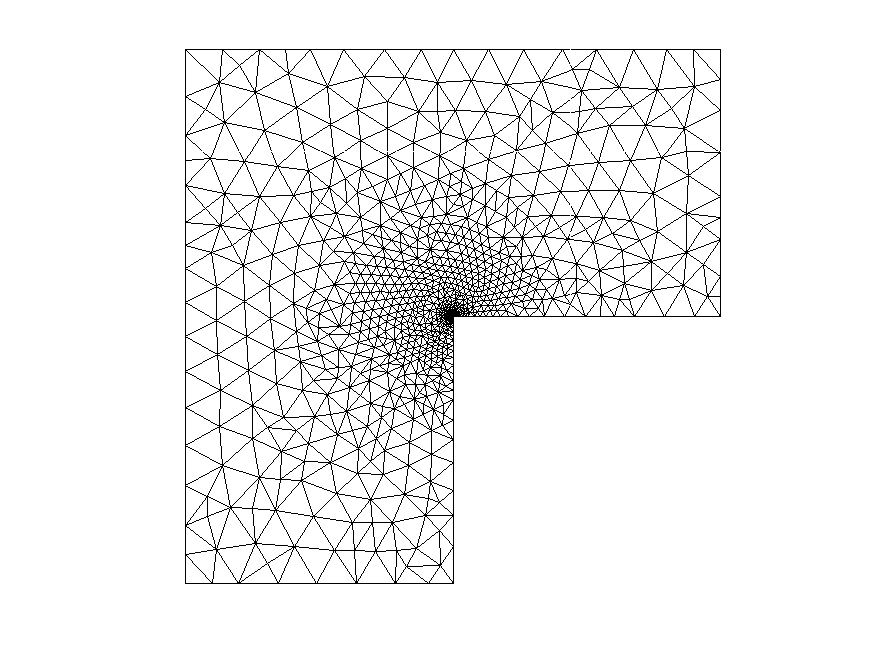}}
		\end{minipage}
		\begin{minipage}{0.3\linewidth}
			\centerline{\includegraphics[width=6cm]{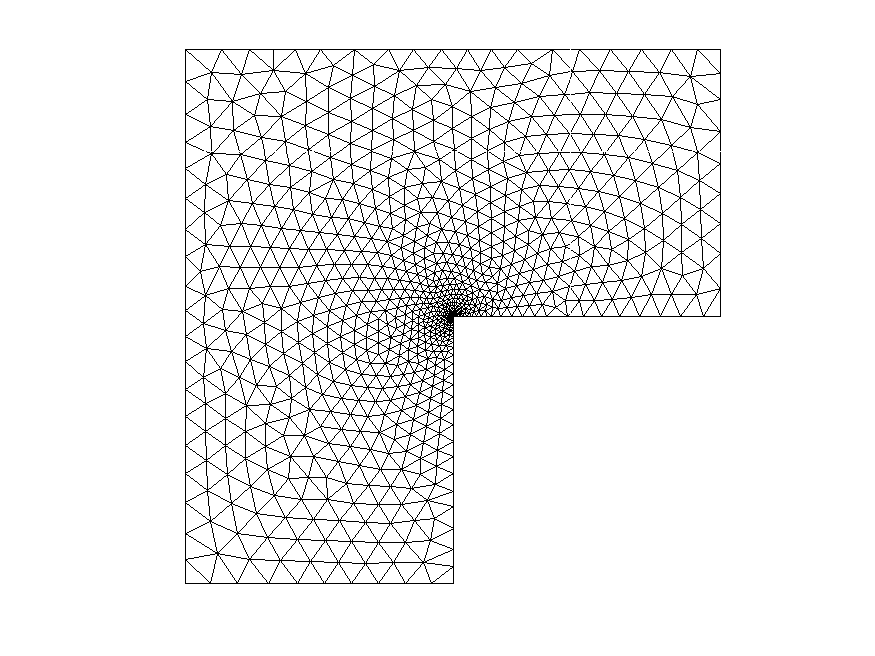}}
		\end{minipage}
		\begin{minipage}{0.3\linewidth}
			\centerline{\includegraphics[width=6cm]{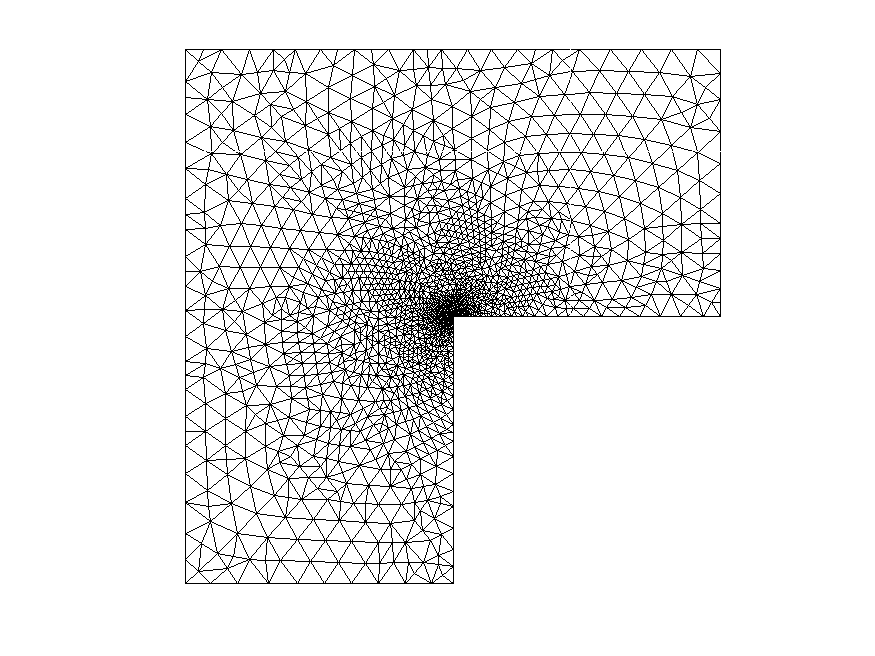}}
		\end{minipage}
		\begin{minipage}{0.3\linewidth}
			\centerline{\includegraphics[width=6cm]{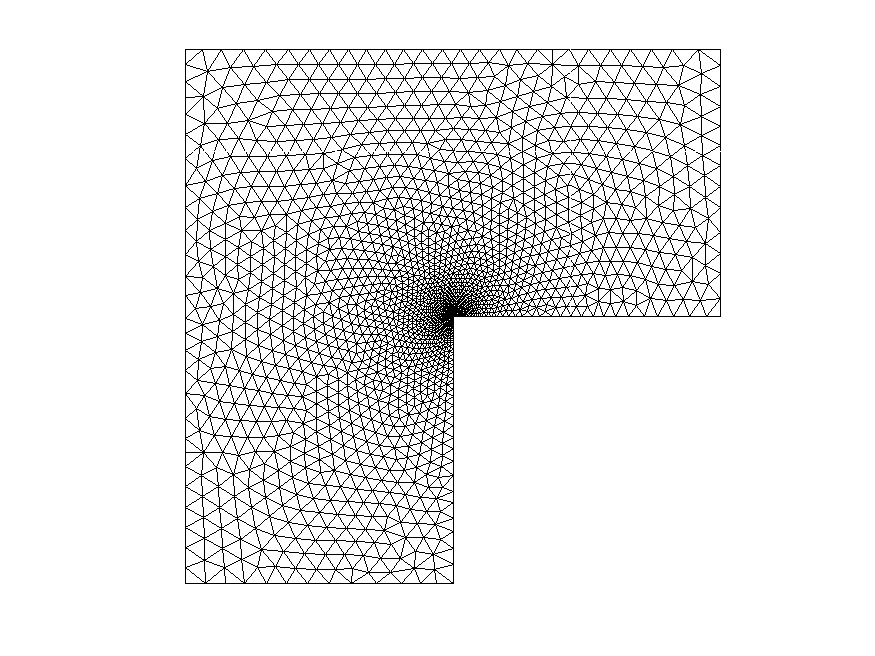}}
		\end{minipage}
		\begin{minipage}{0.3\linewidth}
			\centerline{\includegraphics[width=6cm]{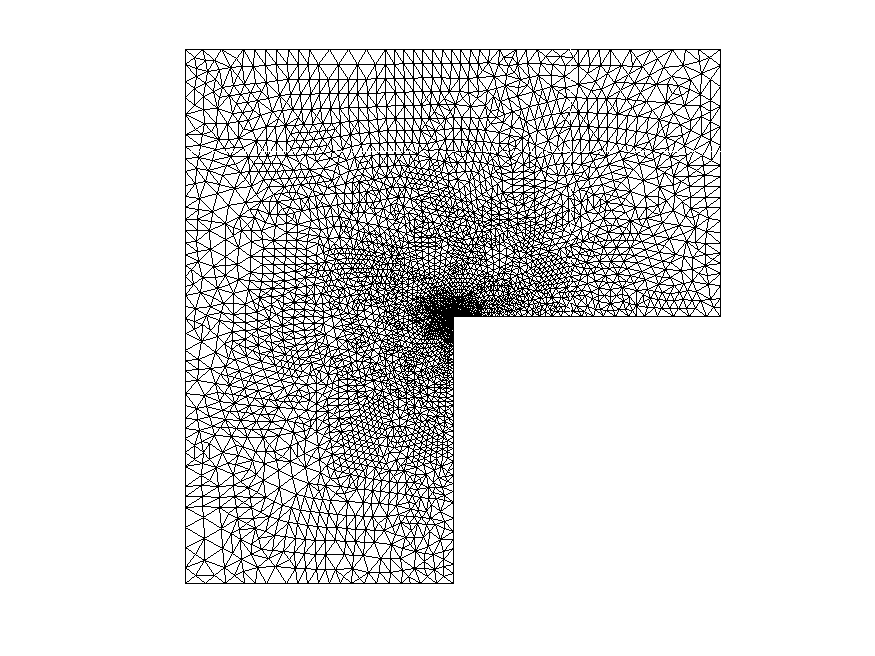}}
		\end{minipage}
		\begin{minipage}{0.3\linewidth}
			\centerline{\includegraphics[width=6cm]{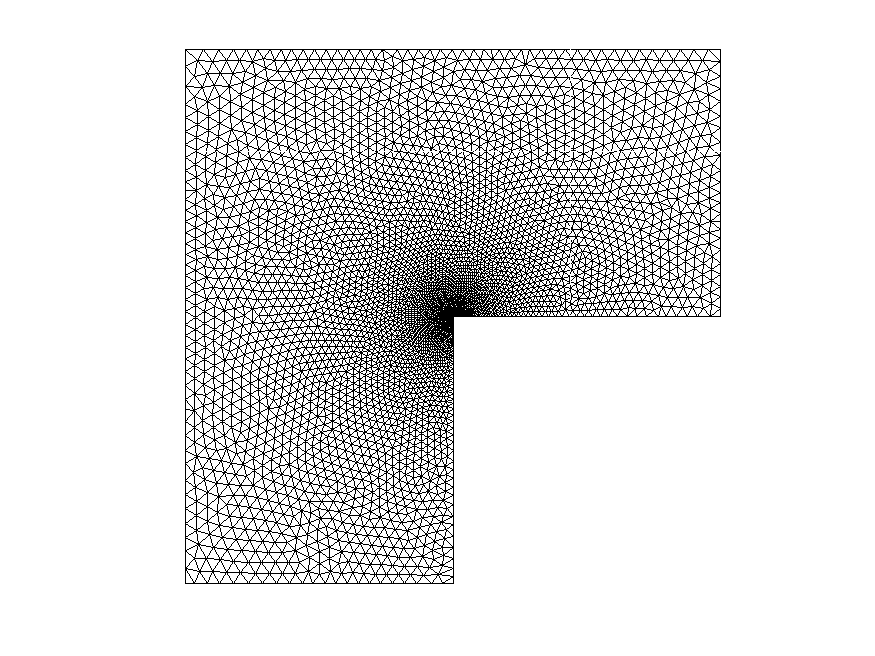}}
		\end{minipage}
		\begin{minipage}{0.3\linewidth}
			\centerline{\includegraphics[width=6cm]{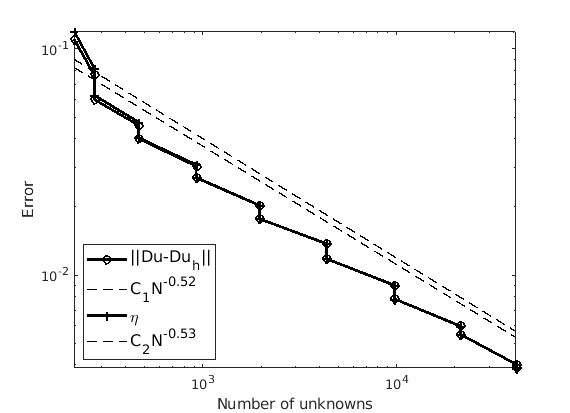}}
		\end{minipage}
		\caption{Example \ref{adap-cvdt}, initial and adaptive meshes, history of gradient error $\|\nabla u-\nabla u_h\|$ and error estimator $\eta_{rec}$.}\label{figLSCT2}
	\end{figure}

	\begin{example}\label{adap-cvdt} 
		We carry out the adaptive finite element method with gradient recovery-based error estimator to the problem presented in Example \ref{eg1}. 
		The adaptive meshes and the history of the numerical error and the gradient recovery type estimator are listed in Figure \ref{figLSCT2}. 
		The initial mesh is generated by a CVDT algorithm with a constant density function, and in each adaptive step we process the loop \textbf{SOLVE} $\rightarrow$ \textbf{ESTIMATE} $\rightarrow$ \textbf{MARK} $\rightarrow$ \textbf{REFINE} $\rightarrow$ \textbf{OPTIMIZATION}. The first four steps are the same in Algorithm \ref{algAFEM}, and in the additional step  \textbf{OPTIMIZATION}, we apply the CVDT algorithm to optimize the mesh with the density function $\rho_h$, which is a piecewise linear function generated from the error estimator. 
		For example, the second mesh shown in Figure \ref{figLSCT2} is a locally refined mesh from the initial mesh, and the third mesh is obtained after step \textbf{OPTIMIZATION}, which takes the second mesh as its input.
		From the error profile, we can find that: i) the error estimator is asymptotically exact; ii) on the adaptive mesh, the finite element approximations on the CVDT mesh are much more accurate than its on the mesh obtained by the newest-vertex bisection refinement.
	\end{example}

	\subsection{Adaptive finite element algorithm}
	Based on the observation shown in Examples \ref{eg1} and \ref{eg3-1}, one should use the asymptotical exact error estimator and as few adaptive iteration steps as possible to develop an efficient AFEM. 
	In the following, we propose a new AFEM, which integrates the finite element high accuracy techniques including the superconvergence and recovery techniques, mesh generation, and optimization in the adaptive procedure. The CVDT mesh optimization is embedded in the mesh adaption to provide high-quality mesh, and then assure the superconvergence property of the recovered gradient and the asymptotical exactness of the error estimator. A tailored adaptive strategy, which could generate high-quality mesh with a target number of vertices, is developed to further improve the efficiency of the adaptive algorithm.
	
	In the CVDT based mesh optimization, one should determine a density function. 
	Following \cite[Sec. 4.1]{Ju2006}, given the error indicators $\{\eta_{h, T}\}_{T\in \mathcal{T}_h}$, we define the piecewise linear (with respect to $\mathcal{T}_h$) density function $\rho_h$ on $\Omega$ such that for any vertex $z_i$ of $\mathcal{T}_h$ 
	\begin{equation}\label{densfun}
		\rho_h(z_i)=\frac{1}{card(\omega_i)}\sum\limits_{T\in\omega_i}\frac{\eta_{h, T}^2}{h_T^4},
	\end{equation}
	where $\omega_i:=\{T\in \mathcal{T}_h| z_i\in \bar{T}\}$.
	
	Given a tolerance $TOL$, we expect the adaptive algorithm to be terminated in less than $7$ steps.
	To achieve this, it is necessary to determine the rate of error reduction in the adaptive procedure. Assume that the error estimator satisfies $\eta_h\approx cN^{-p}$, where the two parameters $c$ and $p$ reflect the rate of error changes with increasing degrees of freedom $N$. 
	We first try five adaptive steps to generate the data $\{(\eta_h^{(i)}, N_i)\}_{i=0}^4$, and obtain the parameters $c$ and $p$ by least square fitting. Then we determine the number of mesh vertices needed to satisfy $cN^{-p}\leq TOL$ and generate the mesh by employing the CVDT method with the density function $\rho_h$ defined by \eqref{densfun}.
	At last, one more adaptive step may be implemented to assure the accuracy of the algorithm. We now summarize our adaptive algorithm in Algorithm \ref{algHATAFEM}.
	
	\normalem
	\begin{algorithm}
		\caption{High accuracy techniques based adaptive finite element algorithm. }\label{algHATAFEM}
		\KwIn{Domain $\Omega$, right hand side function $f$, coefficient matrix $A$, tolerance $TOL$, number of initial mesh vertices $N_0$}
		\KwOut{Sequence of mesh $\{\mathcal{T}_h^{(k)}\}$ and finite element approximations $\{u_h^{(k)}\}$}
		Set $k=0$ and density function $\rho^{(k)}=1$\;
		Generate initial mesh $\mathcal{T}_h^{(k)}$ with $N_{k}$ vertices by Lloyd method with density function $\rho^{(k)}$ \;
		Solve the equation \eqref{lisan} on mesh $\mathcal{T}_h^{(k)}$ to get solution $u_h^{(k)}$ \;
		Calculate local error indicator $\{\eta_{h,T}^{(k)}\}_{T\in\mathcal{T}_h^{(k)}}$ and global error estimator $\eta_{h}^{(k)}$ \;
		\If {$\eta_{h}^{(k)} \leq TOL$}{
			Terminate \;
		}
		\For {$k:=1$ \KwTo $6$}{                
			\If {$k==5$}{
				Calculate the parameters $(c, p)$ in $\eta=cN^{-p}$ by least square fitting with data $\{(\eta_h^{(i)}, N_i)\}_{i=1}^4$ \;
				Determine $N=\left\lceil\sqrt[p]{c/TOL}\right\rceil$ \;
				Set $iteRO=\max\left\{\left\lceil log_2\left(\frac{N}{N_{k-1}}\right)\right\rceil, 1\right\}$ \;
			}
			\Else {
				Set $iteRO=1$ \;
			}
			Construct the density function $\rho^{(k)}$ by \eqref{densfun} \;
			\For {$i:=1$ \KwTo $iteRO$}{   
				Determine $\{\rho(m_i)\}_{i=1}^{N_E^{(k)}}$ and sort them in decreasing order, where $m_i$ denotes the midpoint of edge $E_i$ \;
				Add $\{m_i\}_{i=1}^{n_\theta}$ into the mesh $\mathcal{T}_h^{(k)}$, where 
				\[n_\theta=\max\left\{n| \sum\limits_{i=1}^n\rho(m_i)\leq \frac{1}{2} \sum\limits_{i=1}^{N_E^{(k)}}\rho(m_i)\right\}\]
				and generate the intermediate refined mesh $\tilde{\mathcal{T}}_h^{(k)}$ \;
				Optimize $\tilde{\mathcal{T}}_h^{(k)}$ to obtain $\mathcal{T}_h^{(k)}$ with $N_{k}$ vertices by CfCVDT Algorithm proposed in \cite[ALGORITHM 2]{Ju2006} with density function $\rho^{(k)}$ \;
			}
			Solve the equation \eqref{lisan} on mesh $\mathcal{T}_h^{(k)}$ to get solution $u_h^{(k)}$ \;
			Calculate local error indicator $\{\eta_{h,T}^{(k)}\}_{T\in\mathcal{T}_h^{(k)}}$ and global error estimator $\eta_{h}^{(k)}$ \;
			\If {$\eta_{h}^{(k)}\leq TOL$ \textbf{or} $k==6$}{
				Terminate \;
			}
		}
	\end{algorithm}

	\section{Numerical results}\label{sec:num}
	This section reports three numerical experiments that exhibit a variety of types of singularities to verify the efficiency and robustness of the proposed adaptive algorithm.  
	
\begin{figure}[ht]
		\begin{minipage}{0.3\linewidth}
			\centerline{	\includegraphics[width=5.5cm]{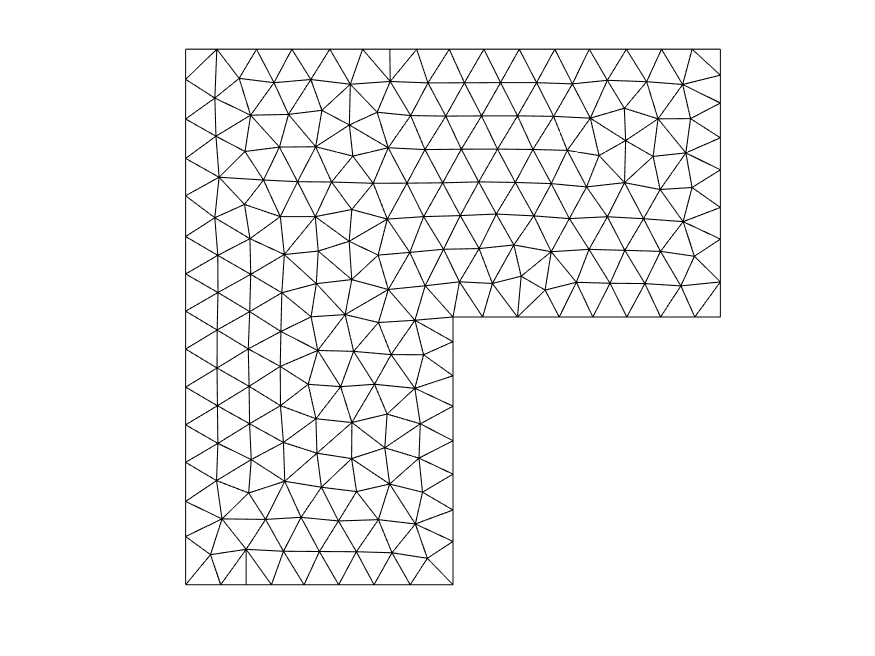}}
			\centerline{(a)}
		\end{minipage}
		\begin{minipage}{0.3\linewidth}
			\centerline{	\includegraphics[width=5.5cm]{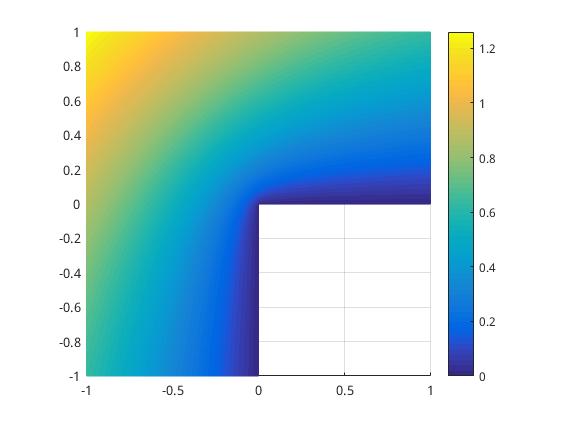}}
			\centerline{(b)}
		\end{minipage}
		\begin{minipage}{0.3\linewidth}
			\centerline{	\includegraphics[width=5.5cm]{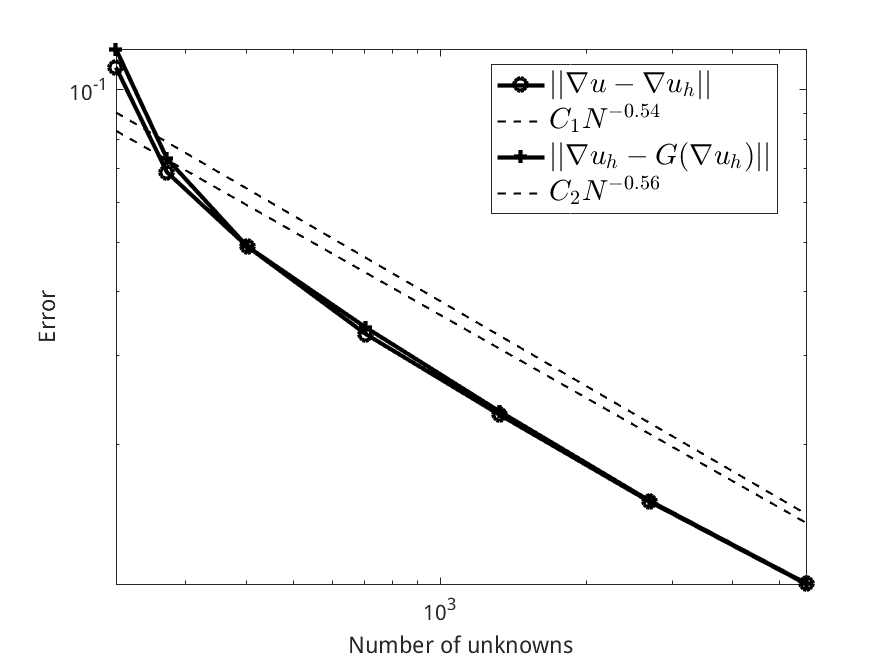}}
			\centerline{(c)}
		\end{minipage}
		\caption{Example \ref{ne1}, (a)  initial  mesh; (b) numerical solution; (c)
			history of the error and estimator.}\label{figure5}
	\end{figure}

\begin{figure}[!htbp]
		\begin{minipage}{0.3\linewidth}
			\centerline{	\includegraphics[width=5.5cm]{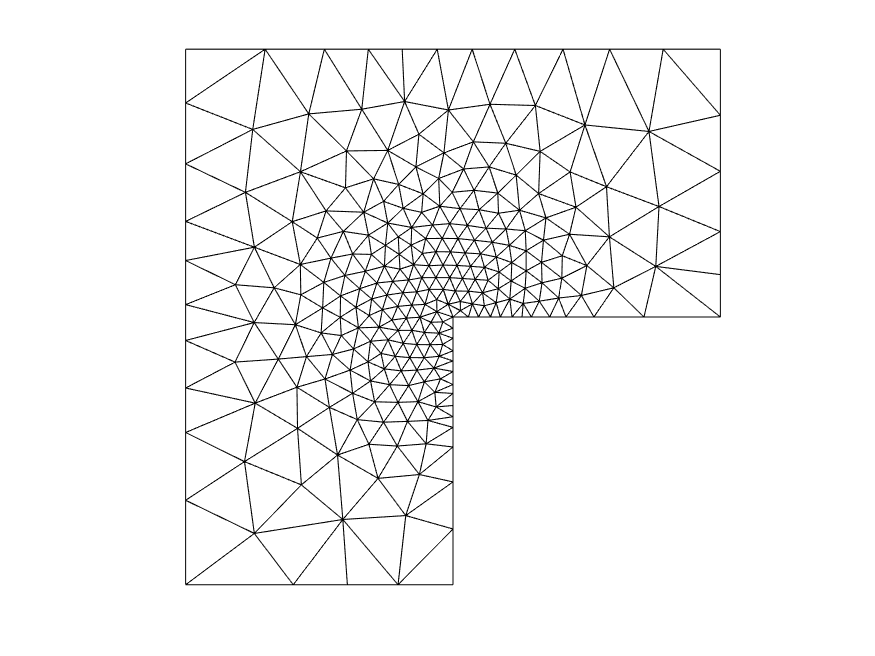}}
		\end{minipage}
		\begin{minipage}{0.3\linewidth}
			\centerline{	\includegraphics[width=5.5cm]{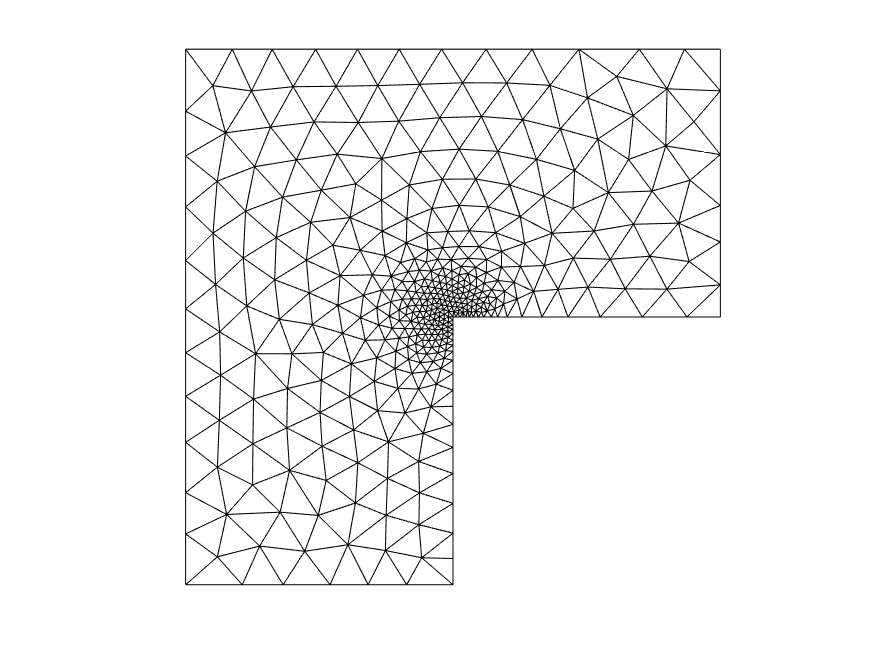}}
		\end{minipage}\begin{minipage}{0.3\linewidth}
			\centerline{	\includegraphics[width=5.5cm]{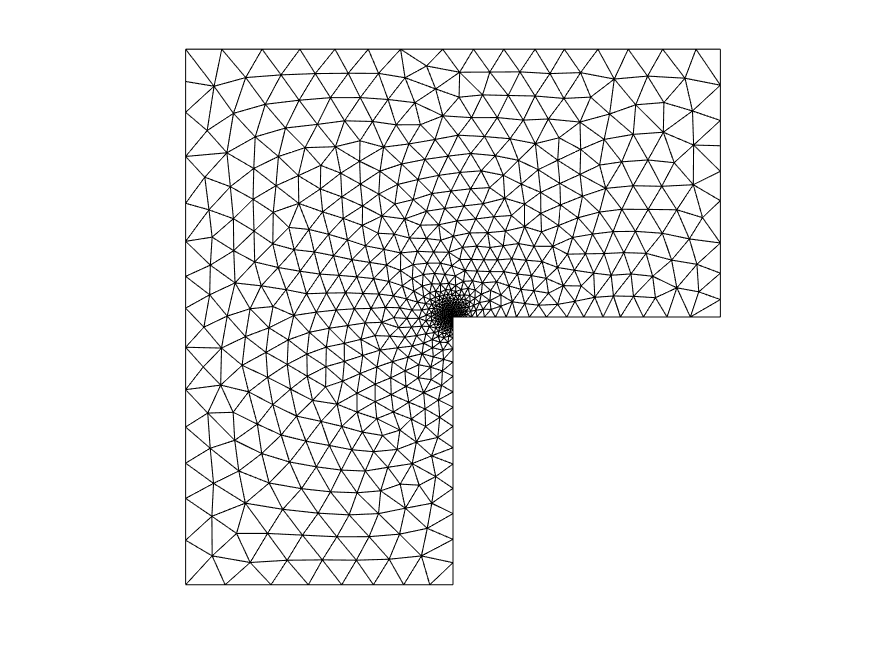}}
		\end{minipage}	
		\begin{minipage}{0.3\linewidth}
			\centerline{	\includegraphics[width=5.5cm]{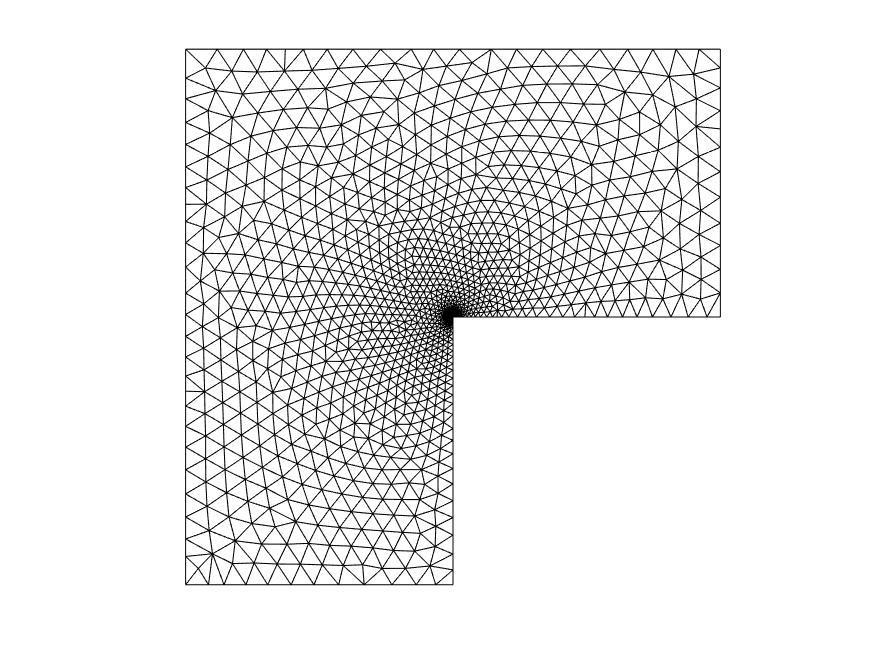}}
		\end{minipage}
		\begin{minipage}{0.3\linewidth}
			\centerline{	\includegraphics[width=5.5cm]{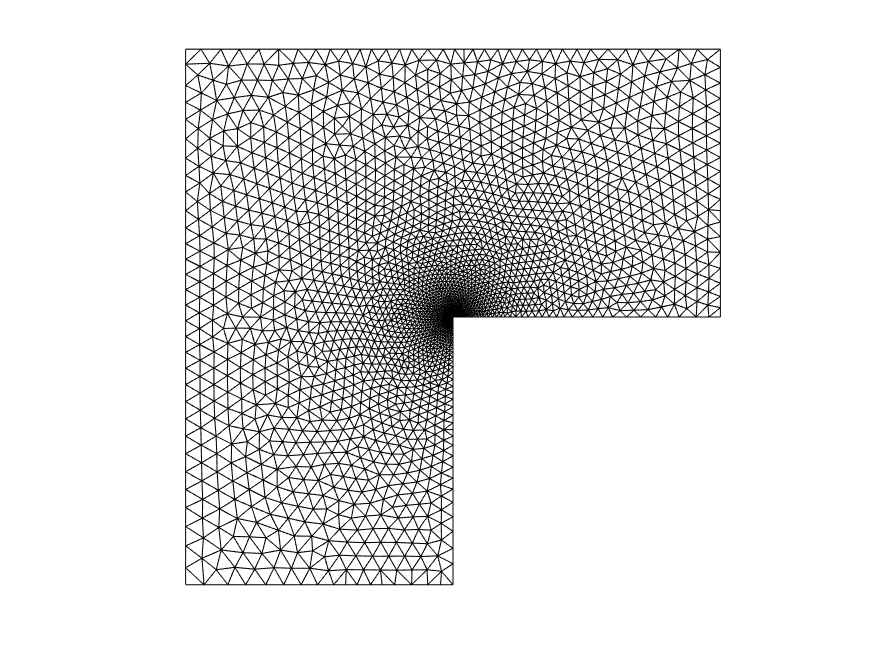}}
		\end{minipage}
		\begin{minipage}{0.3\linewidth}
			\centerline{	\includegraphics[width=5.5cm]{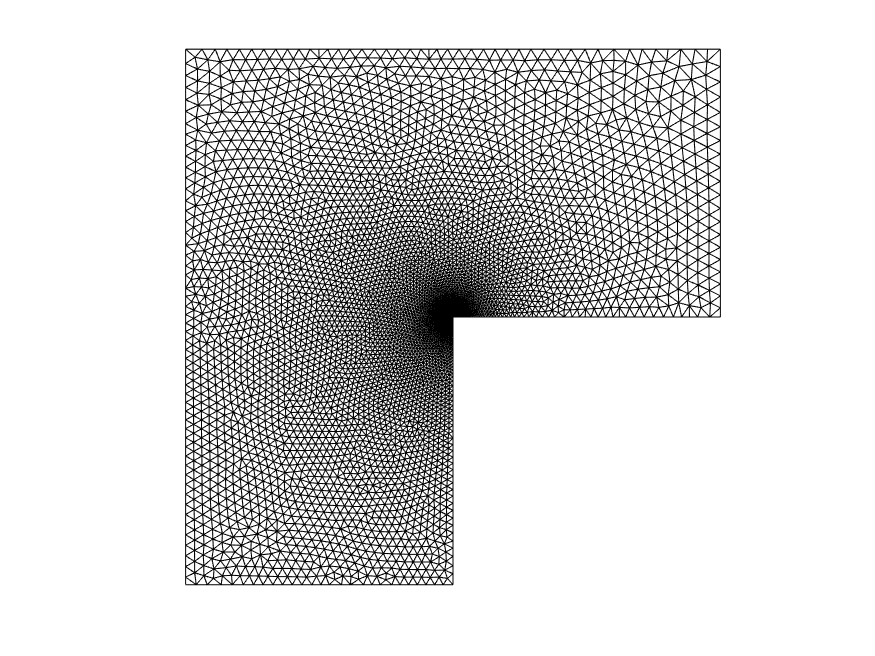}}
		\end{minipage}
		\caption{Example \ref{ne1}, adaptive meshes.}\label{fig_eg1_mesh}
	\end{figure}

	\begin{example}\label{ne1}[Corner singularity problem] We apply Algorithm \ref{algHATAFEM} with $TOL=0.01$ to the L-shape domain problem reported in Example \ref{eg1}. The initial mesh is shown in Figure \ref{figure5} (a), and the corresponding adaptive meshes are shown in Figure \ref{fig_eg1_mesh}. It shows that the error indicator successfully guides the mesh refinement near the reentrant corner $(0, 0)$, and all the meshes generated in the adaptive procedure are of high quality. From Figure \ref{figure5} (c), we can see that the error and the estimator achieve the optimal convergence, and the estimator is asymptotically exact. 
		We list the data of mesh vertices, error $\|\nabla u-\nabla u_h\|$, and the recovery type error estimator $\eta_{rec}$ in Table \ref{LSCVcan}. Compared to the data listed in Table \ref{Lshpaedata} of Algorithm \ref{algAFEM} with the residual type estimator, Algorithm \ref{algHATAFEM} requires only $7$ iterations and stops at the mesh with $5671$ vertices, thus saving significantly on computation cost.
	\end{example}
	
	\begin{table}[ht]
		\caption{Example \ref{ne1}, data of the mesh vertices, error and recovery type error estimator.}\label{LSCVcan} \centering
		\begin{tabular}{cccc} \hline
			$k$  &     $N$      & $||\nabla u-\nabla u_h||$ & $\eta_{rec}$ \\\hline
			1    &     216      & 1.1053e-01  & 1.2012e-01 \\
			2    &     275      & 6.8665e-02  & 7.3026e-02\\
			3    &     402      & 4.9057e-02  & 4.9062e-02 \\
			4    &     705      & 3.2947e-02  & 3.3949e-02 \\
			5    &     1326     & 2.2869e-02  & 2.3215e-02 \\
			6    &     2702     & 1.5404e-02  & 1.5467e-02 \\
			7    &     5671     & 1.0601e-02  & 1.0629e-02 \\\hline
		\end{tabular}
	\end{table}
	
\begin{figure}[!ht]
		\begin{minipage}{0.3\linewidth}
			\centerline{	\includegraphics[width=5.5cm]{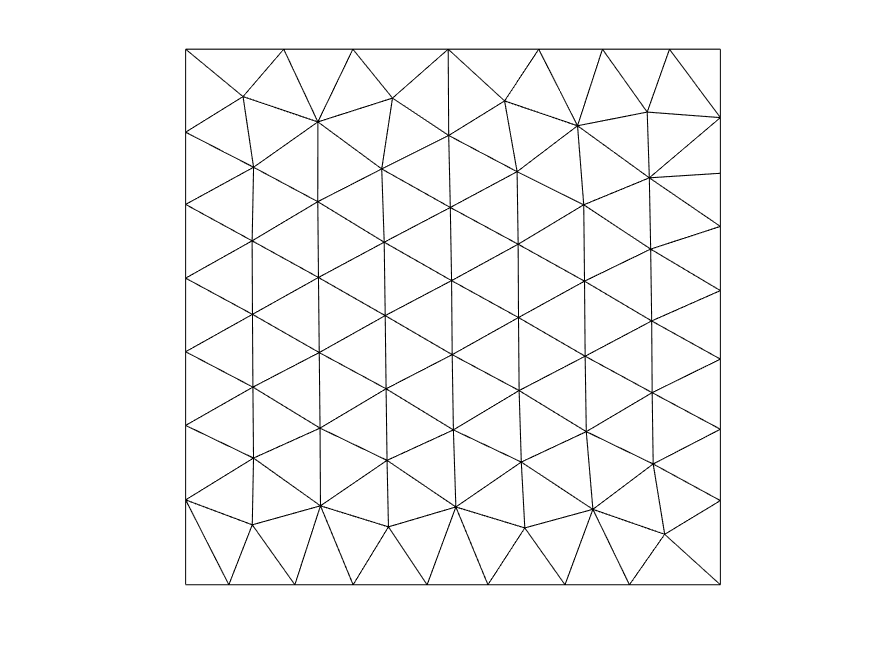}}
			\centerline{(a)}
		\end{minipage}
		\begin{minipage}{0.3\linewidth}
			\centerline{	\includegraphics[width=5.5cm]{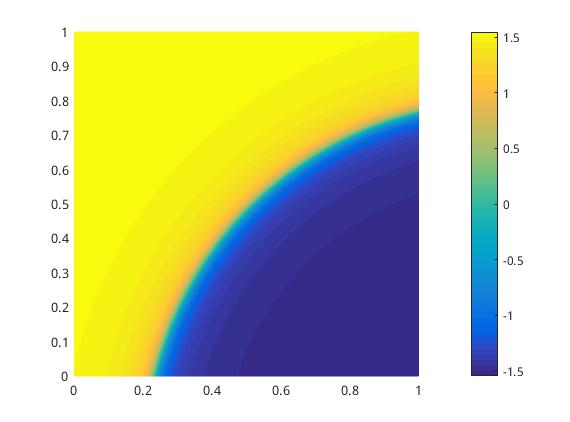}}
			\centerline{(b)}
		\end{minipage}
		\begin{minipage}{0.3\linewidth}
			\centerline{	\includegraphics[width=5.5cm]{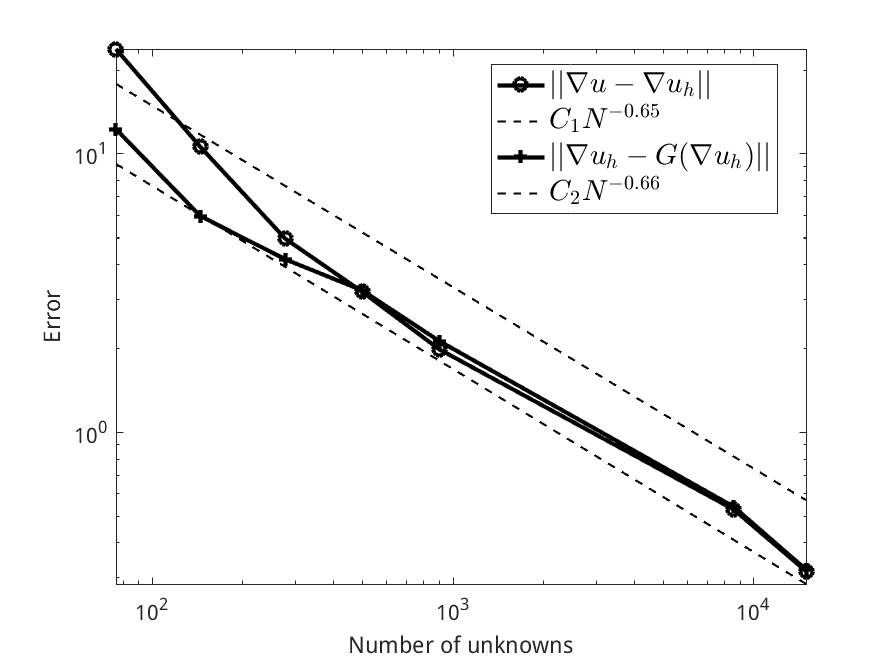}}
			\centerline{(c)}
		\end{minipage}
		\caption{Example \ref{eg_inner}, (a) initial mesh; (b) numerical solution; (c)
			history of the error and estimator.}\label{figLSCWW}
	\end{figure}

\begin{figure}[!ht]
		\begin{minipage}{0.3\linewidth}
			\centerline{	\includegraphics[width=5.5cm]{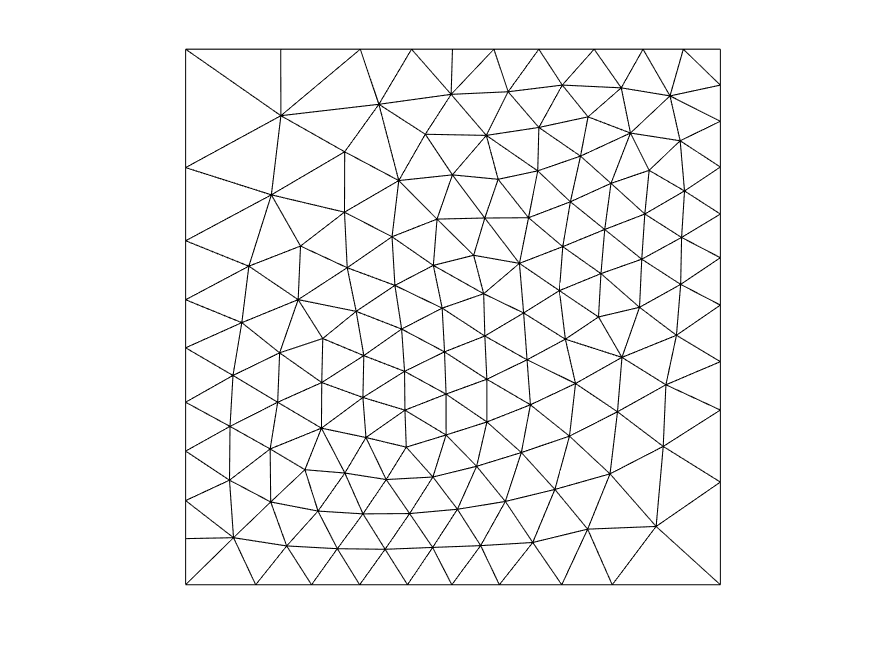}}
		\end{minipage}
		\begin{minipage}{0.3\linewidth}
			\centerline{	\includegraphics[width=5.5cm]{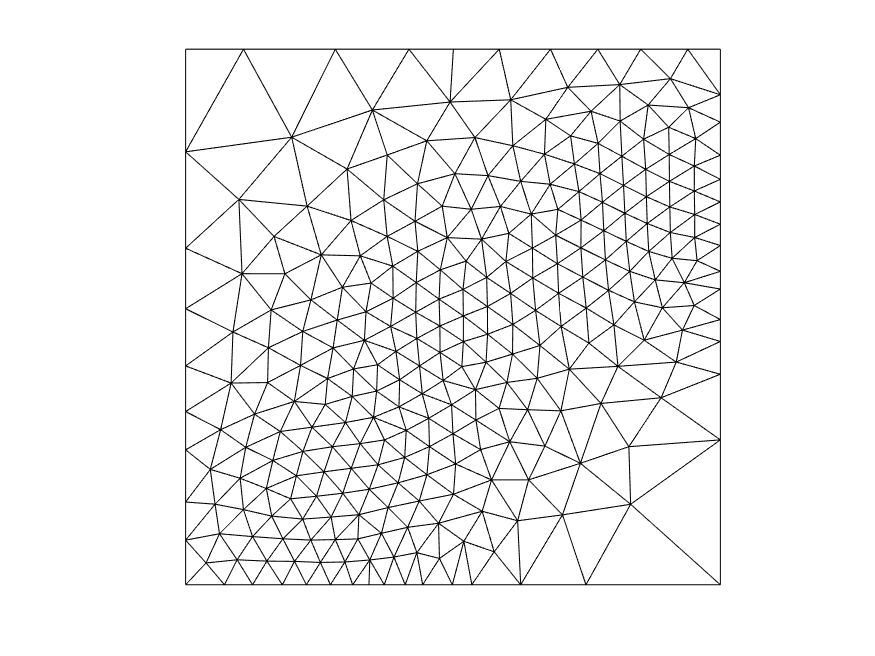}}
		\end{minipage}\begin{minipage}{0.3\linewidth}
			\centerline{	\includegraphics[width=5.5cm]{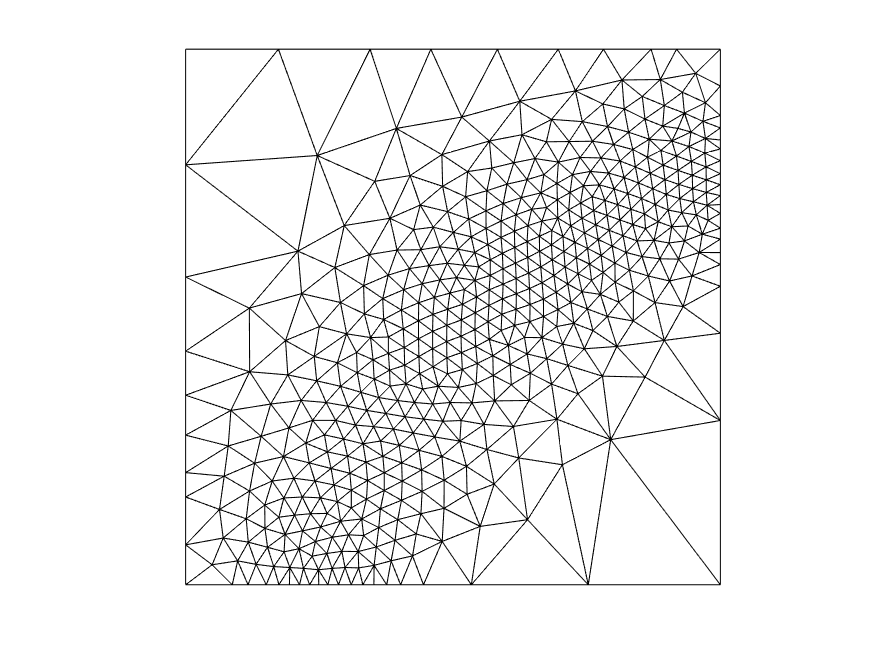}}
		\end{minipage}
		
		\begin{minipage}{0.3\linewidth}
			\centerline{	\includegraphics[width=5.5cm]{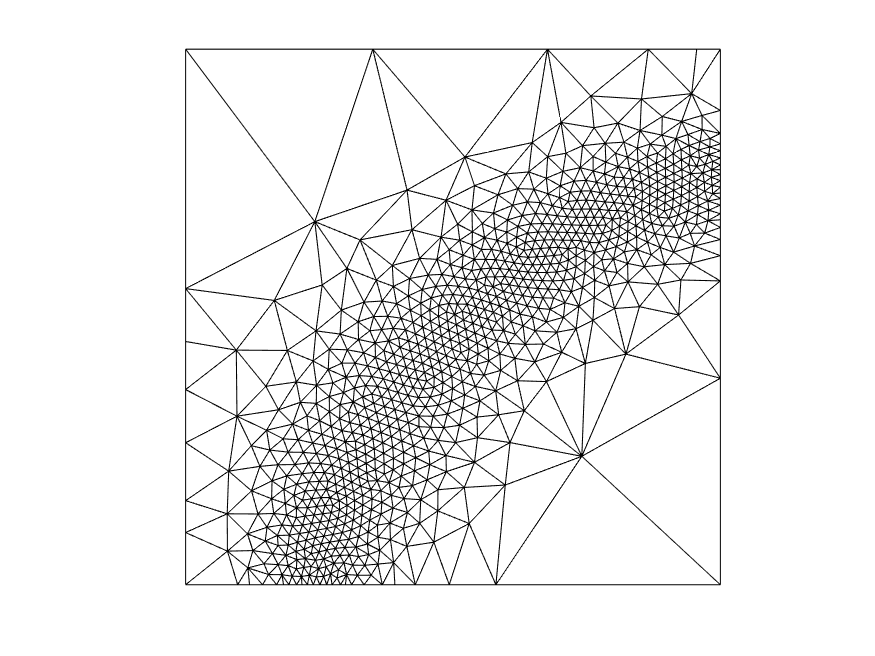}}
		\end{minipage}
		\begin{minipage}{0.3\linewidth}
			\centerline{	\includegraphics[width=5.5cm]{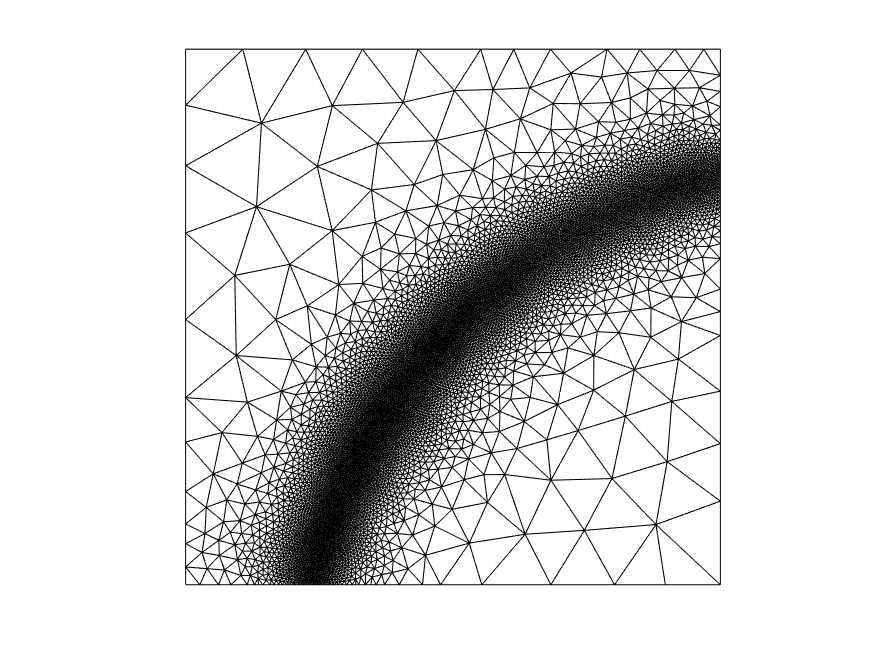}}
		\end{minipage}
		\begin{minipage}{0.3\linewidth}
			\centerline{	\includegraphics[width=5.5cm]{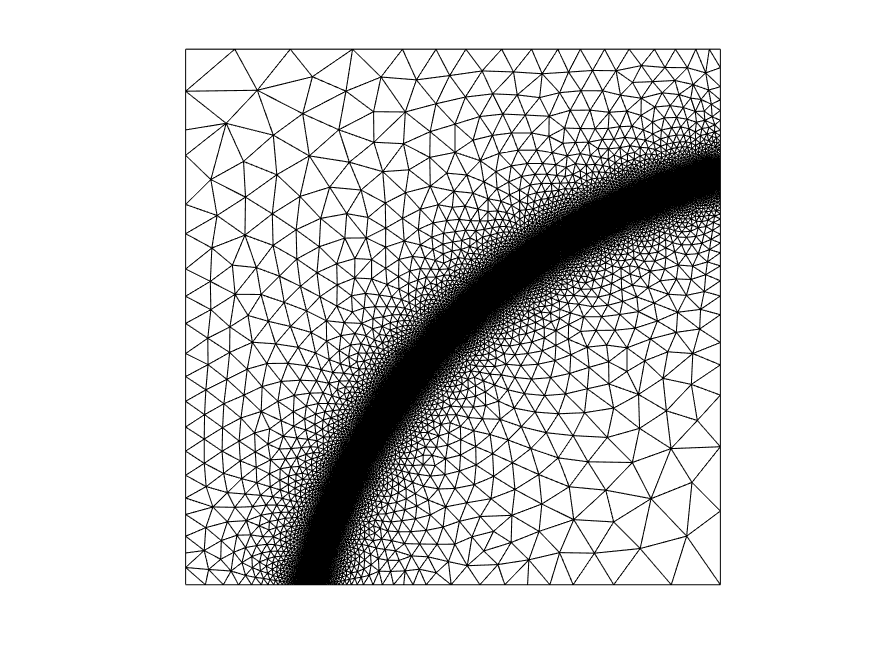}}
		\end{minipage}
		\caption{Example \ref{eg_inner}, adaptive meshes.}\label{figLSCTT}
	\end{figure}

	\begin{example}[Inner layer problem]\label{eg_inner}
		Let $\Omega = [0,1]^2$, we consider the
		Poisson equation with Dirichlet boundary condition given by a smooth solution 
		\[u = atan(S(\sqrt{(x-1.25)^2+(y+0.25)^2}-\pi/3)),\]
		where $S=60$ reflects the steepness of the inner slope and the source term $f$ is obtained from $u$.
		
		Apply Algorithm \ref{algHATAFEM} with $TOL=0.5$ and an initial uniform CVDT mesh with $76$ vertices. Figure \ref{figLSCWW} displays the initial mesh, numerical solution, and convergence history of error and estimator, respectively. Figure \ref{figLSCWW} (c) demonstrates the decays of error $\|\nabla u-\nabla u_h\|$ and error estimator $\eta_{rec}$ are optimal, the rates are approximate $O(N^{-0.65})$, and the recovery type error estimator is asymptotically exact. The adaptive meshes shown in Figure \ref{figLSCTT} demonstrate the capability of the error estimator to track the high gradient domain. We can see that the mesh refinement arises along with the inner layer which indicates that the singularities are perfectly captured by the adaptive algorithm. 
		Table \ref{LSCVcan} lists the data of mesh vertices, error $\|\nabla u-\nabla u_h\|$ and the recovery type error estimator $\eta_{rec}$. We can see that: i) Algorithm \ref{algHATAFEM} requires only $7$ iterations and stopped at the mesh with $14960$ vertices, ii) the error estimator is asymptotically exact, iii) the vertices of mesh at the adaptive steps $5$ and $6$ are $904$ and $8593$, respectively, this shows that our new Algorithm \ref{algHATAFEM} can generate a mesh with target number of vertices at once. 
	\end{example}

	\begin{table}[!ht]
		\caption{Example \ref{eg_inner}, data of the mesh vertices, error and recovery type error estimator.}\label{LSCVcann} \centering
		\begin{tabular}{cccc} \hline
			$k$ &     $N$      & $||\nabla u-\nabla u_h||$ & $\eta_{rec}$\\ \hline
			1   &     76       & 2.3688e+01  & 1.2216e+01 \\
			2   &     145      & 1.0542e+01  & 5.9498e+00 \\
			3   &     278      & 4.9312e+00  & 4.1612e+00 \\
			4   &     500      & 3.1838e+00  & 3.2206e+00 \\
			5   &     904      & 1.9760e+00  & 2.1153e+00 \\
			6   &     8593     & 5.2594e-01  & 5.3796e-01 \\
			7   &     14960    & 3.1564e-01  & 3.1983e-01 \\ \hline
		\end{tabular}
	\end{table}

\begin{figure}[!ht]
		\begin{minipage}{0.3\linewidth}
			\centerline{	\includegraphics[width=5.5cm]{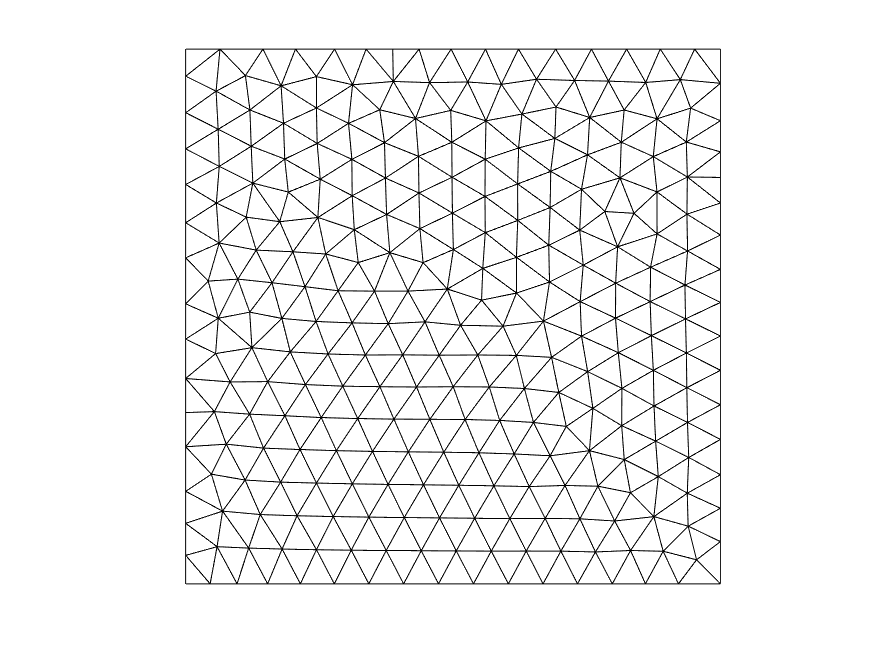}}
			\centerline{(a)}
		\end{minipage}
		\begin{minipage}{0.3\linewidth}
			\centerline{	\includegraphics[width=5.5cm]{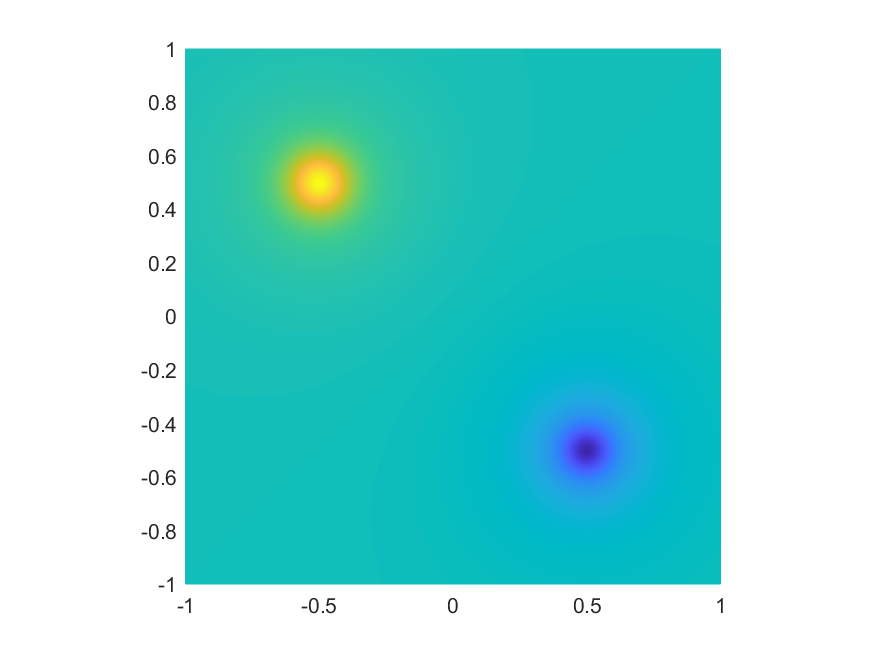}}
			\centerline{(b)}
		\end{minipage}
		\begin{minipage}{0.3\linewidth}
			\centerline{	\includegraphics[width=5.5cm]{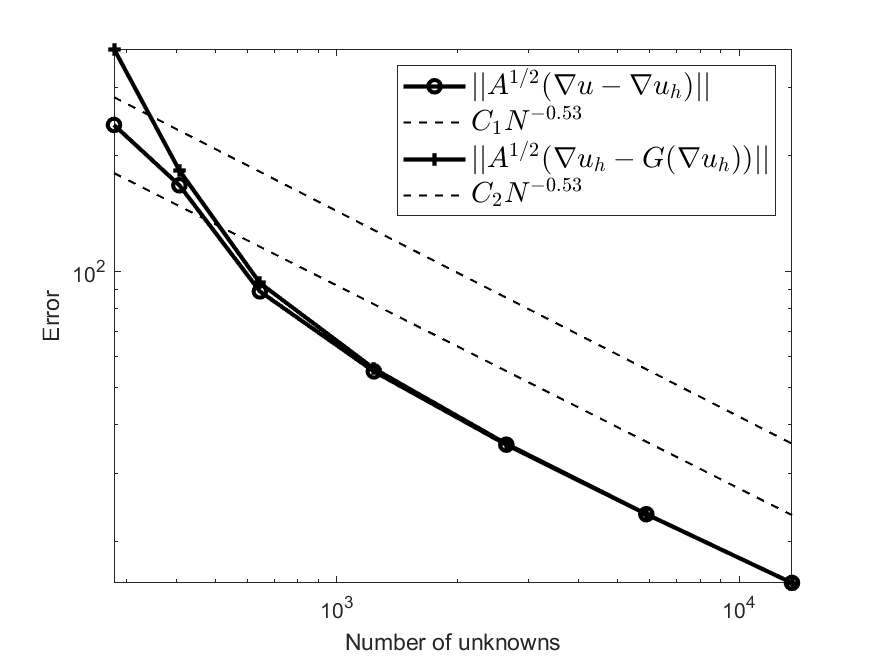}}
			\centerline{(c)}
		\end{minipage}
		\caption{Example \ref{eg_sharp}, (a) initial mesh; (b) numerical solution; (c)
			history of the error and estimator.}\label{figSPWW}
	\end{figure}

\begin{figure}[!ht]
		\begin{minipage}{0.3\linewidth}
			\centerline{	\includegraphics[width=5.5cm]{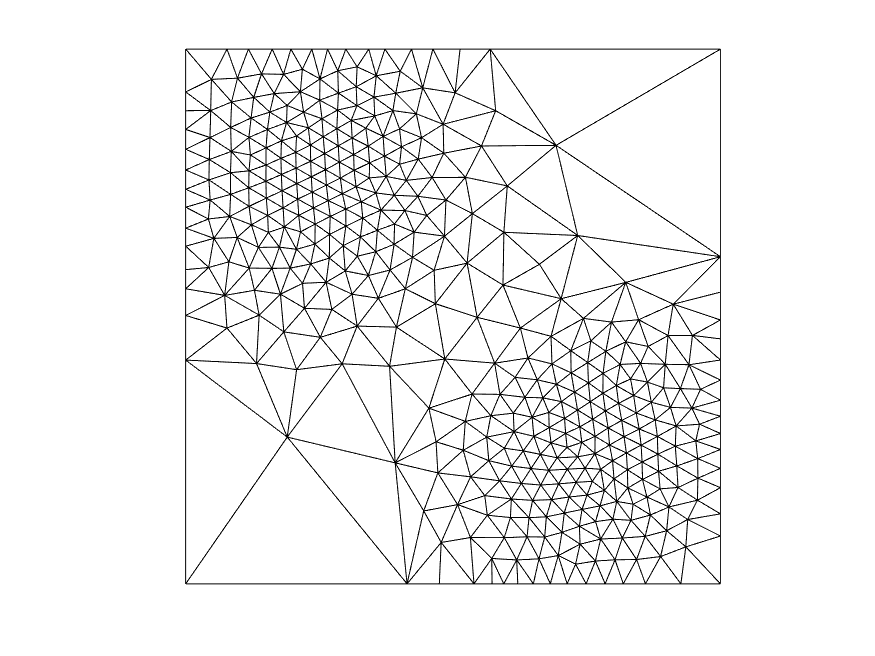}}
		\end{minipage}
		\begin{minipage}{0.3\linewidth}
			\centerline{	\includegraphics[width=5.5cm]{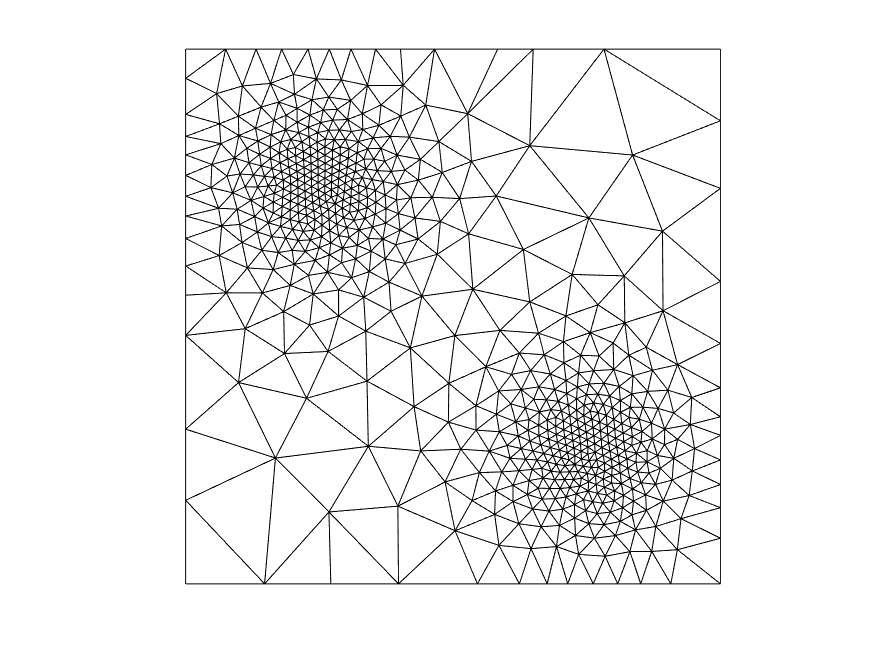}}
		\end{minipage}\begin{minipage}{0.3\linewidth}
			\centerline{	\includegraphics[width=5.5cm]{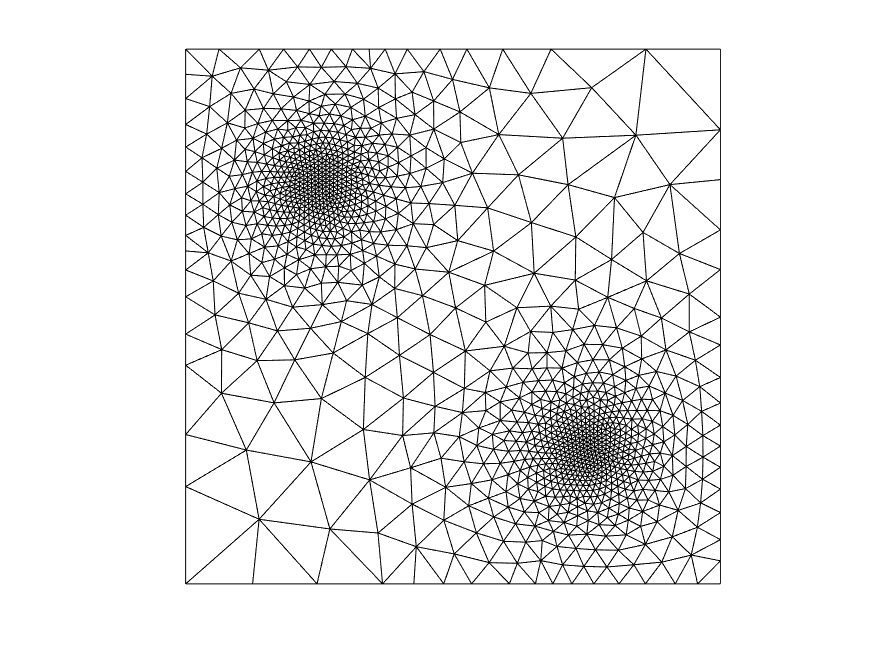}}
		\end{minipage}
		
		\begin{minipage}{0.3\linewidth}
			\centerline{	\includegraphics[width=5.5cm]{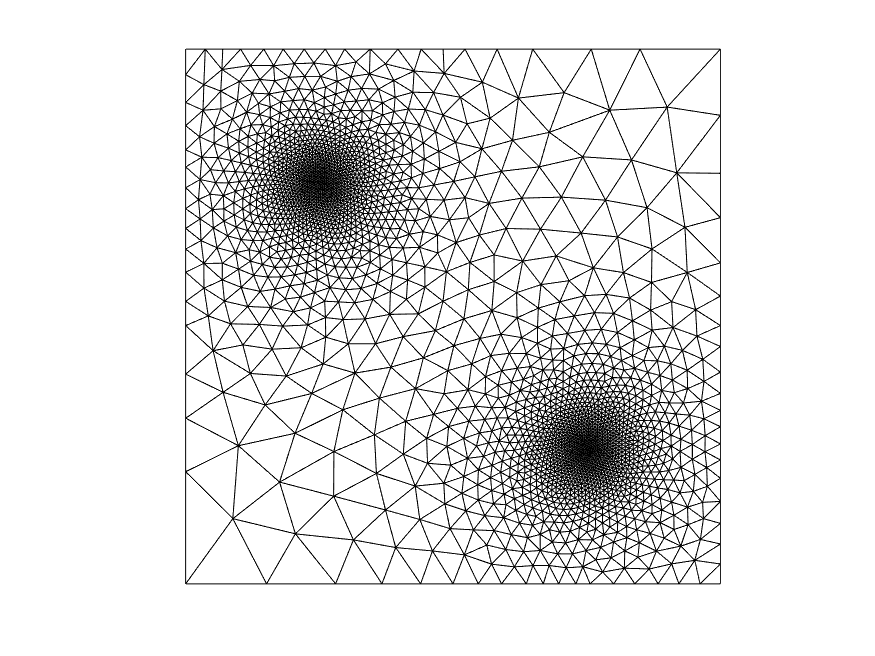}}
		\end{minipage}
		\begin{minipage}{0.3\linewidth}
			\centerline{	\includegraphics[width=5.5cm]{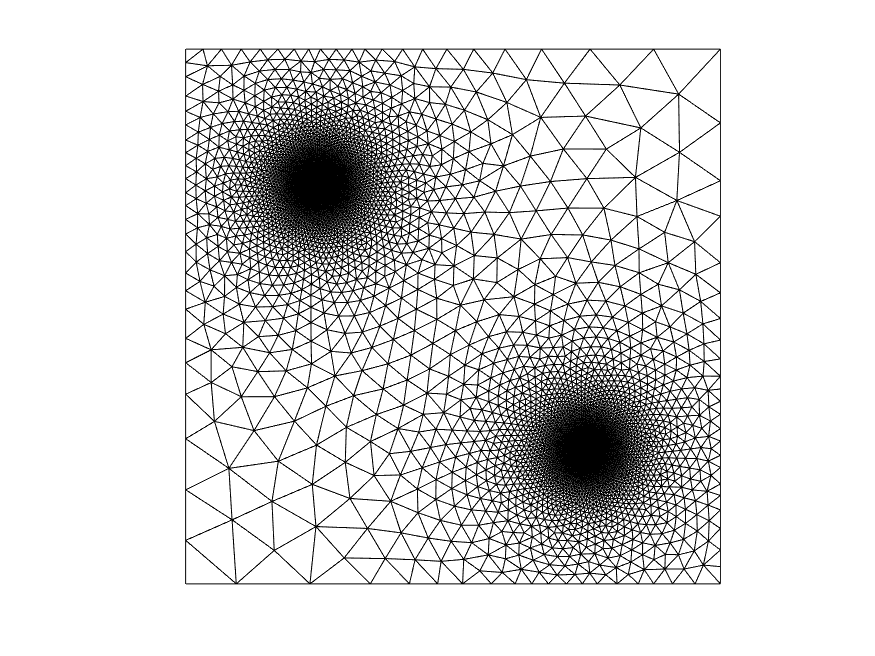}}
		\end{minipage}
		\begin{minipage}{0.3\linewidth}
			\centerline{	\includegraphics[width=5.5cm]{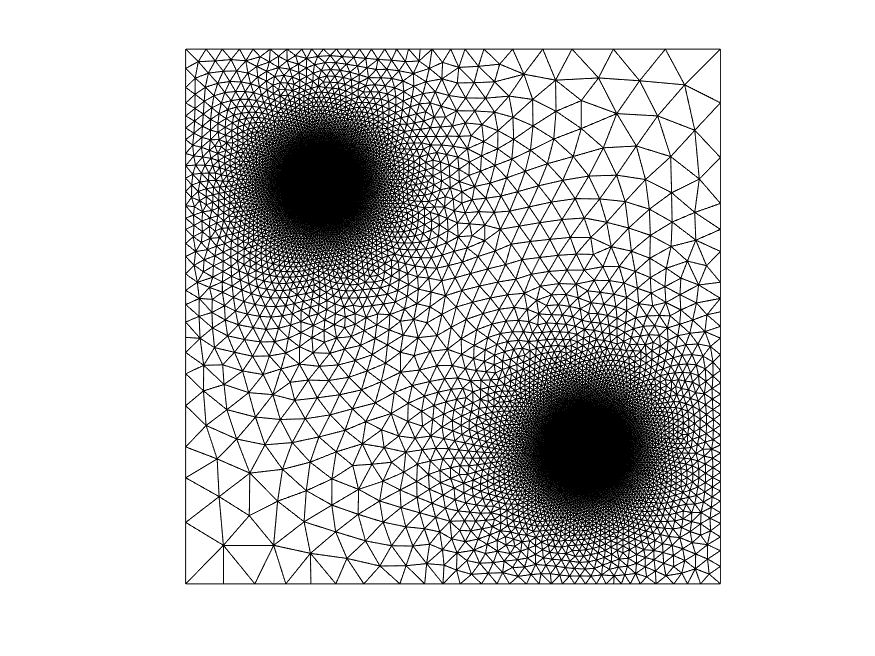}}
		\end{minipage}
		\caption{Example \ref{eg_sharp}, adaptive meshes.}\label{figSPTT}
	\end{figure}

	\begin{example}[Peak problem]\label{eg_sharp}
		Let $\Omega = [-1,1]^2$, we consider the problem 
		\[\left\{\begin{aligned}
			-\nabla\cdot(A\nabla u) &=f\quad \text{in}~\Omega, \\
			u&=g\quad \text{on}~\partial \Omega,
		\end{aligned}\right.\]
		with continuous diffusion coefficient $A = 10\cos y$. The exact solution is
		$$u=\frac{1}{(x+0.5)^2+(y-0.5)^2+0.01}-\frac{1}{(x-0.5)^2+(y+0.5)^2+0.01},$$
		where $f$ and $g$ if obtained from $u$. 
		
		The exact solution $u$ is a smooth function that achieves its maximum value $99\frac{101}{201}$ at the point $(-0.5,0.5)$ and its minimum value $-99\frac{101}{201}$ at $(0.5,-0.5)$, but decays quickly away from its extrema and thus has large gradients near these two points. 
		Noting that the equation contains the coefficients $A$, we modify the error estimator accordingly as
		\[\eta_{T, rec}^2=\|A^{1/2}(G(\nabla u_h)-\nabla u_h)\|_{0, T},\quad \eta_{rec}=\|A^{1/2}(G(\nabla u_h)-\nabla u_h)\|.\]
		We start with an initial uniform CVDT mesh with $280$ vertices and $TOL=20$. 
		Figure \ref{figSPWW} plots the initial mesh, numerical solution, and the errors, respectively, we can see that the error and the estimator arrive at the optimal convergence rate, and the estimator is asymptotically exact. 
		Toward the two peaks 
		at $(-0.5, 0.5)$ and $(0.5, -0.5)$, adaptive refined meshes are shown in Figure \ref{figSPTT}. We see clearly that the meshes generated by the adaptive algorithm are all of high quality.
		We list the data of mesh vertices and errors in Table \ref{figSPCC}. The adaptive process is terminated in $7$ steps. 
		Both the errors and the mesh refinement demonstrate that the adaptive algorithm leads to a very effective convergence procedure.
	\end{example}
	
	\begin{table}[!ht]
		\caption{Example \ref{eg_sharp}, data of the mesh vertices, error and recovery type error estimator.}\label{figSPCC} \centering
		\begin{tabular}{cccc} \hline
			$k$    &     $N$    & $||A^{\frac{1}{2}}(\nabla u-\nabla u_h)||$ & $\eta_{rec}$                   \\ \hline
			1      &     280    & 2.3968e+02  & 3.7703e+02 \\
			2      &     407    & 1.6718e+02  & 1.8260e+02 \\
			3      &     645    & 8.8735e+01  & 9.3325e+01 \\
			4      &     1238   & 5.5050e+01  & 5.5968e+01 \\
			5      &     2641   & 3.5545e+01  & 3.5720e+01 \\
			6      &     5884   & 2.3469e+01  & 2.3490e+01 \\
			7      &     13537  & 1.5585e+01  & 1.5576e+01 \\ \hline
		\end{tabular}
	\end{table}
	
	
	\section*{Acknowledgments}
	Yi's research was partially supported by the NSFC Project (12431014), 
	Project of Scientiﬁc Research Fund of the Hunan Provincial Science and Technology Department (2024ZL5017), 
	and Program for Science and Technology Innovative Research Team in Higher Educational Institutions of Hunan Province. Liu's research was supported by the NSFC Project (12301473) and the Natural Science Special Project of the Shaanxi Provincial Department of Education (23JK0564).
	

\begin{thebibliography}{10}
		
		\bibitem{AD2000}
		M.~ Ainsworth and J.T.~ Oden.
		\newblock A posteriori error estimation in finite element analysis.
		\newblock {\em Wiley Interscience, New York},  2000.
		
		
		\bibitem{AD2000}
		D.~Arnold, A.~Mukherjee, and L.~Pouly.
		\newblock Locally adapted tetrahedral meshes using bisection.
		\newblock {\em SIAM J. Sci. Comput.}, 22(2):431--448, 2000.
		
		\bibitem{babuska2012}
		I.~Babu{\v{s}}ka, J.~Flaherty, W.~Henshaw, J.~Hopcroft, J.~Oliger, and
		T.~Tezduyar.
		\newblock {\em Modeling, mesh generation, and adaptive numerical methods for
			partial differential equations}, volume~75.
		\newblock Springer Science \& Business Media, 2012.
		
		\bibitem{br}
		I.~Babu{\v{s}}ka and W.~Rheinboldt.
		\newblock A-posteriori error estimates for the finite element method.
		\newblock {\em Int. J. Numer. Meth. Eng.}, 12:1597--1615, 1978.
		
		\bibitem{Babu_1978}
		I.~Babu{\v{s}}ka and W.~Rheinboldt.
		\newblock Error estimates for adaptive finite element computations.
		\newblock {\em SIAM J. Numer. Anal.}, 15(4):736--754, 1978.
		
		\bibitem{bank1998}
		R.~Bank.
		\newblock {\em PLTMG: A software package for solving elliptic partial
			differential Equations: Users' Guide 8.0}.
		\newblock SIAM, 1998.
		
		\bibitem{Bernardi2000}
		C.~Bernardi and R.~Verfürth.
		\newblock Adaptive finite element methods for elliptic equations with
		non-smooth coefficients.
		\newblock {\em Numer. Math.}, 85:579--608, 2000.
		
		\bibitem{BCNV2024}
		A.~Bonito, C.~Canuto, R.~Nochetto, and A.~Veeser.
		\newblock Adaptive finite element methods.
		\newblock {\em Acta Numer.}, 33:163 -- 485, 2024.
		
		\bibitem{cz}
		Z.~Cai and S.~Zhang.
		\newblock Recovery-based error estimator for interface problems: conforming
		linear elements.
		\newblock {\em SIAM J. Numer. Anal.}, 47:2132--2156, 2009.
		
		\bibitem{chen2004}
		L.~Chen and J.~Xu.
		\newblock Optimal delaunay triangulations.
		\newblock {\em J. Comput. Math.}, pages 299--308, 2004.
		
		\bibitem{cn}
		Z.~Chen and R.H. Nochetto.
		\newblock Residual type a posteriori error estimates for elliptic obstacle
		problems.
		\newblock {\em Numer. Math.}, 84:527--548, 2000.
		
		\bibitem{D1996}
		W.~D{\"o}rfler.
		\newblock A convergent adaptive algorithm for {P}oisson’s equation.
		\newblock {\em SIAM J. Numer. Anal.}, 33(3):1106--1124, 1996.
		
		\bibitem{yan2001}
		L.~Du and N.~Yan.
		\newblock Gradient recovery type a posteriori error estimate for finite element
		approximation on non-uniform meshes.
		\newblock {\em Adv. Comput. Math.}, 14:175--193, 2001.
		
		\bibitem{Du1999}
		Q.~Du, V.~Faber, and M.~Gunzburger.
		\newblock Centroidal voronoi tessellations: Applications and algorithms.
		\newblock {\em SIAM Rev.}, 41(4):637--676, 1999.
		
		\bibitem{du2002}
		Q.~Du and M.~Gunzburger.
		\newblock Grid generation and optimization based on centroidal {V}oronoi
		tessellations.
		\newblock {\em Appl. Math. Comput.}, 133(2-3):591--607, 2002.
		
		\bibitem{fa}
		F.~Fierro and A.~Veeser.
		\newblock A posteriori error estimators, gradient recovery by averaging, and
		superconvergence.
		\newblock {\em Numer. Math.}, 103:267--298, 2006.
		
		\bibitem{hy}
		Y.~Huang, K.~Jiang, and N.~Yi.
		\newblock Some weighted averaging methods for gradient recovery.
		\newblock {\em Adv. Appl. Math. Mech.}, 4:131--155, 2012.
		
		\bibitem{HQW2008}
		Y.~Huang, H.~Qin, and D.~Wang.
		\newblock Centroidal {Voronoi} tessellation‐based finite element
		superconvergence.
		\newblock {\em Int. J. Numer. Method. Engrg.}, 76:1819–1839, 2008.
		
		\bibitem{Huang2011}
		Y.~Huang, H.~Qin, D.~Wang, and Q.~Du.
		\newblock Convergent adaptive finite element method based on centroidal
		{V}oronoi tessellations and superconvergence.
		\newblock {\em Commun. Comput. Phys.}, 10:339--370, 2011.
		
		\bibitem{HWY2017}
		Y.~Huang, L.~Wang, and N.~Yi.
		\newblock Mesh quality and more detailed error estimate of finite element
		method.
		\newblock {\em Numer. Math. Theor. Meth. Appl.}, 10(2):420--436, 2017.
		
		\bibitem{hy1}
		Y.~Huang and N.~Yi.
		\newblock The superconvergent cluster recovery method.
		\newblock {\em J. Sci. Comput.}, 44:301--322, 2010.
		
		\bibitem{joe1989}
		B.~Joe.
		\newblock Three-dimensional triangulations from local transformations.
		\newblock {\em SIAM J. Sci. Stat. Comput.}, 10(4):718--741, 1989.
		
		\bibitem{Ju2006}
		L.~Ju, M.~Gunzburger, and W.~Zhao.
		\newblock Adaptive finite element methods for elliptic {PDE}s based on
		conforming centroidal {V}oronoi–{D}elaunay triangulations.
		\newblock {\em SIAM J. Sci. Comput.}, 28:2023--2053, 2006.
		
		\bibitem{KOKO2015650}
		Jonas Koko.
		\newblock A {Matlab} mesh generator for the two-dimensional finite element
		method.
		\newblock {\em Appl. Math. Comput.}, 250:650--664, 2015.
		
		\bibitem{LXYC2024}
		Y.~Liu, J.~Xiao, N.~Yi, and H.~Cao.
		\newblock Gradient recovery-based a posteriori error estimator and adaptive
		finite element method for elliptic equations.
		\newblock {\em submited}, 2024.
		
		\bibitem{DJ1990}
		D.~Mavriplis.
		\newblock Adaptive mesh generation for viscous flows using triangulation.
		\newblock {\em J. Comput. Phys.}, 90(2):271--291, 1990.
		
		\bibitem{MKN2005}
		K.~Mekchay and R.H. Nochetto.
		\newblock Convergence of adaptive finite element methods for general second
		order linear elliptic {PDE}s.
		\newblock {\em SIAM J. Numer. Anal.}, 43:1803--1827, 2005.
		
		\bibitem{mn}
		P.~Morin, R.H. Nochetto, and K.~Siebert.
		\newblock Convergence of adaptive finite element methods.
		\newblock {\em SIAM Rev.}, 44:631--658, 2002.
		
		\bibitem{PJ1992}
		J.~Peraire, J.~Peiro, and K.~Morgan.
		\newblock Adaptive re-meshing for three-dimensional compressible flow
		computations.
		\newblock {\em J. Comput. Phys.}, 103(2):269--285, 1992.
		
		\bibitem{rivara1984}
		M.~Rivara.
		\newblock Mesh refinement processes based on the generalized bisection of
		simplices.
		\newblock {\em SIAM J. Numer. Anal.}, 21(3):604--613, 1984.
		
		\bibitem{Karel2010}
		Karel Segeth.
		\newblock A review of some a posteriori error estimates for adaptive finite
		element methods.
		\newblock {\em Math. Comput. Simulat}, 80:1589--1600, 2010.
		
		\bibitem{DW2005}
		D.~Wang and Q.~Du.
		\newblock Mesh optimization based on centroidal voronoi tessellation.
		\newblock {\em Int. J. Numer. Anal. Model.}, 2:100–114, 2005.
		
		\bibitem{yz}
		N.~Yan and A.~Zhou.
		\newblock Gradient recovery type a posteriori error estimates for finite
		element approximations on irregular meshes.
		\newblock {\em Comput. Method. Appl. M.}, 190:4289--4299, 2001.
		
		\bibitem{za}
		Z.~Zhang and A.~Naga.
		\newblock A new finite element gradient recovery method: superconvergence
		property.
		\newblock {\em SIAM J. Sci. Comput.}, 26:1192--1213, 2005.
		
		\bibitem{zz}
		O.~Zienkiewicz and J.~Zhu.
		\newblock A simple error estimator and adaptive procedure for practical
		engineering analysis.
		\newblock {\em Int. J. Numer. Meth. Eng.}, 24:337--357, 1987.
		
	\end{thebibliography}

\end{document}